\documentclass{article}
\usepackage[utf8]{inputenc}
\usepackage{amsmath}
\usepackage{amssymb}
\usepackage{multirow}
\usepackage{algorithm}
\usepackage{algorithmic}
\usepackage{caption}
\usepackage{subcaption}
\usepackage{float}
\usepackage{graphicx}
\usepackage[english]{babel}
\usepackage[utf8]{inputenc}
\usepackage{adjustbox}
\usepackage{url}
\usepackage{outlines}

\title{Reduced order model of a convection-diffusion equation using Proper Orthogonal Decomposition}
\author{Neelakantan Padmanabhan \\
npadman@g.clemson.edu}
\date{}

\begin{document}

\maketitle
\section{Abstract}
In this work, a numerical simulation of 1D Burgers' equation is developed using finite difference method and a reduced order model (ROM) of the simulation is developed using proper orthogonal decomposition (POD). The objective of this work is to provide an introduction of the POD method to researchers interested in computational fluid dynamics (CFD). This work discusses a physical interpretation of the POD method, its strengths and shortcomings and an implementation of the algorithm that may be extended to 2D, 3D Burgers' equation and other non-linear partial differential equations (PDE) of this class, to develop models for more complex systems.

\section{Introduction}
Burgers' equation is a fundamental nonlinear partial differential equation that finds applications in the areas of thermo-fluids and CFD \cite{burgers1,burgers2,burgers3}. It is a generalized convection-diffusion equation and falls under the same category of equations as the Navier-Stokes equations. Numerical solution and reduced order modeling of Burgers' equation, is of great interest in the CFD community. Numerical solution of Burgers' equation and other PDEs of this class can be computationally very expensive \cite{dissertation}. Thus, reduced order models are highly desirable since they can capture the essence of the phenomenon with a fraction of the computational resources. POD / PCA (principal component analysis) is a well established approach to develop reduced order models \cite{POD2,POD3,POD5,POD6,POD7,POD8,POD9,POD10,PCA1}. For a given data, POD extracts the modes that contain the most dominant characteristics of the data. In a function space, these modes form the orthonormal basis that can be used to reconstruct the data. Ideally, POD modes can be constructed from any simulation or experimental data and if the governing equations of the system are known, Galerkin projection may be applied on a subset of this basis to create a ROM, to predict the time evolution of the system. While this approach results in accurate ROM for systems with high diffusivity, the predictions for systems with low diffusivity especially for long simulation times are observed to be less accurate. A number of improvements to this method have been proposed, that help in developing more accurate models. Some of these include a large eddy simulation like approach \cite{dissertation}, where the large scale effects are modeled by POD and the small scale effects are modeled by eddy viscosity type model \cite{POD1}, goal-oriented POD \cite{POD4}, and discrete empirical interpolation method (DEIM) \cite{POD_DEIM1,POD_DEIM2}. \\
\section{Numerical Simulation of Burgers' equation}
1D Burgers' equation of the form presented in Eq.\ref{eq_burgers} is solved using finite difference approach. 
\begin{equation}
\frac{\partial U}{\partial t} = - U \frac{\partial U}{\partial x} + \nu \frac{\partial ^2 U}{\partial x^2} + Q,
\label{eq_burgers}
\end{equation}
where, $U=U(x,t)$ represents a field variable, $t$ represents time, $x$ represents the spatial vector, $\nu$ represents the diffusivity, and $Q$ a source term. The generic Burgers' equation incorporates a non-linear convection term, and a linear diffusion term. An additional source term is considered in this work. A second order central finite difference scheme is used for spatial discretization of the diffusion term, while a second order upwind scheme is employed for the nonlinear convection term. An explicit forward Euler scheme is used for time stepping. 
\begin{equation}
\begin{aligned}
    \frac{U(x,t+\Delta t) - U(x,t)}{\Delta t} = -U(x,t) \frac{3 U(x,t)-4 U(x-\Delta x,t) + U(x-2\Delta x,t)}{2\Delta x} \\ + \nu \frac{U(x+\Delta x,t) - 2 U(x,t) + U(x-\Delta x,t)}{\Delta x^2} + Q
\end{aligned}
\label{FD_disc}
\end{equation}
The equations are solved on an equally spaced Cartesian grid. Constant diffusivity values are used in the simulations. CFL conditions for velocity and diffusivity are used to determine the minimum time step. Periodic boundary conditions are used in the domain. The following initial conditions are used for the field variable and source term, 
\begin{equation}
\begin{aligned}
    U(x,0) = U_0 sin(x),\\
    Q(x,0) = Q_0 sin(x),
\end{aligned}
\end{equation}
where, $u_0=0.01, Q_0=0.1$. Simulations are run for various configurations of number of grid points and duration ($t_{final}$). The results from the numerical simulations are presented in Figs. \ref{fig_sim_case1}-\ref{fig_sim_case5}.
\section{Reduced Order Modeling}
An overview of a few concepts of linear algebra that are central to the POD method is presented in Appendix \ref{linalgterm}. A brief review of this section might be useful in developing a better understanding of the method. 
\subsection{A physical interpretation of POD}
\begin{figure}
    \centering
	\begin{subfigure}{0.45\linewidth}
		\includegraphics[width=\textwidth]{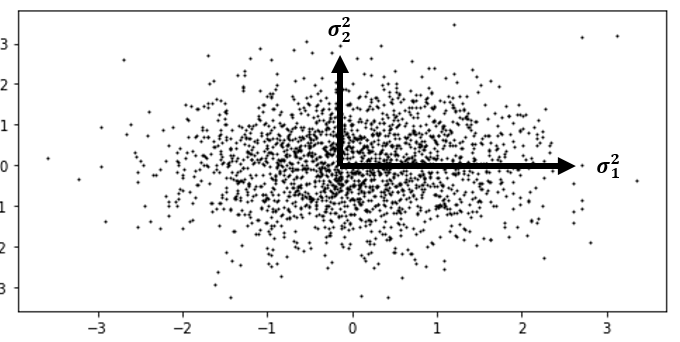}
	\end{subfigure} \hspace{1em}
	\caption{A simulated plot of randomly scattered data. $\sigma_1 ^2, \sigma_2 ^2$ represent the directions of highest and lowest variance.}
	\label{fig:gen}
\end{figure}

Solutions of Burgers' equation ($U(x,t)=U \in \mathbb{R}^{m \times n}$, with $m$ spatial points and $n$ temporal points), can be expressed analytically in terms of the basis vectors $x$ and $t$, in a physical space. While it is desirable, it is not always convenient or viable to obtain an analytical solution in this space. Alternatively, this solution can be also expressed in terms of a different set of basis vectors in a function space. This alternative basis can be determined by singular value decomposition (SVD) of the solution matrix $U$ or by eigen decomposition of the covariance matrices of the solution matrix ($Cov_1(U)=UU^*$ and $Cov_2(U)=U^*U$). The covariance matrices, represent the variance exhibited by the elements of the solution matrix $U$ and the covariance between the pairs of dataset. It is useful in separating structured relationships in a matrix of random variables. When eigen decomposition of the covariance matrix is performed, it results in a factorization of the form $Cov(U)=V\Lambda V^{-1}$. Here, $V$ represents the eigenvectors or the principle components and $\Lambda$ represents the eigenvalues of the covariance matrix. Geometrically, this can be interpreted as a linear transformation represented by $Cov(U)$ which when applied to the vectors $V$, only results in $V$ scaling up or down by a factor of $\Lambda$. Since $Cov(U)$ represent the variance in the data, the vectors $V$ indicate the directions along which the variance in the data is the highest or the lowest and $\Lambda$ represents the magnitude of significance of the given eigenvector. In other terms, the variance indicates how the energy or information of the solution is distributed and the eigenvectors indicate the directions in which the distribution is significant. A graphical illustration of this is presented in Fig.\ref{fig:gen}. In a function space, the vectors $V$ form an orthonormal basis. In many cases (except the cases with very low diffusivity in the Burgers' equation), it is observed that the first few basis vectors (or modes) typically contain most of the energy of the solution. Hence, by projecting $U$ on these first few modes, the highest energy solutions can be recovered. By Galerkin projection of a weak form of the governing equation onto a sub-space of this basis, a reduced order model can be obtained, from which the time evolution of the system can be predicted.

\subsection{Formal description of POD}
\begin{figure}[!htb]
	\centering
	\begin{subfigure}{0.4\linewidth}
		\includegraphics[width=\textwidth]{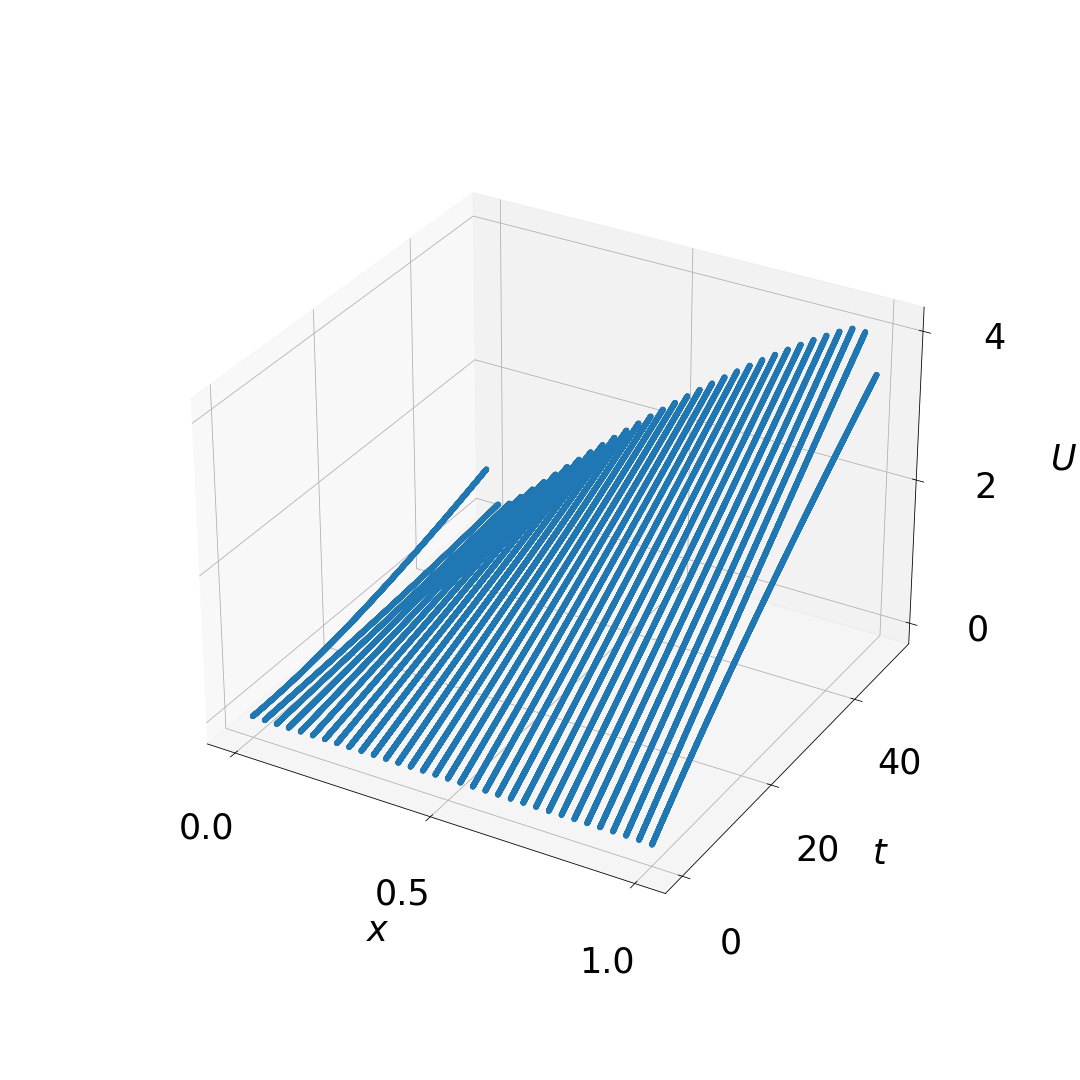} 
		\caption{Simulation with $t_{final}=50$, $N_x=32$, $CFL=0.2$, $\nu=0.00001$.}
		\label{fig_sim_case5}
	\end{subfigure} \hspace{1em}
	\begin{subfigure}{0.4\linewidth}
		\includegraphics[width=\textwidth]{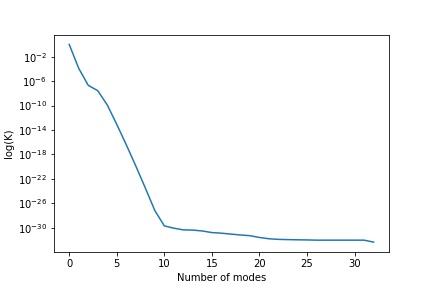} 
		\caption{Energy of dominant modes in the simulation}
		\label{fig_energy_case5}
	\end{subfigure} \\
	\begin{subfigure}{0.4\linewidth}
		\includegraphics[width=\textwidth]{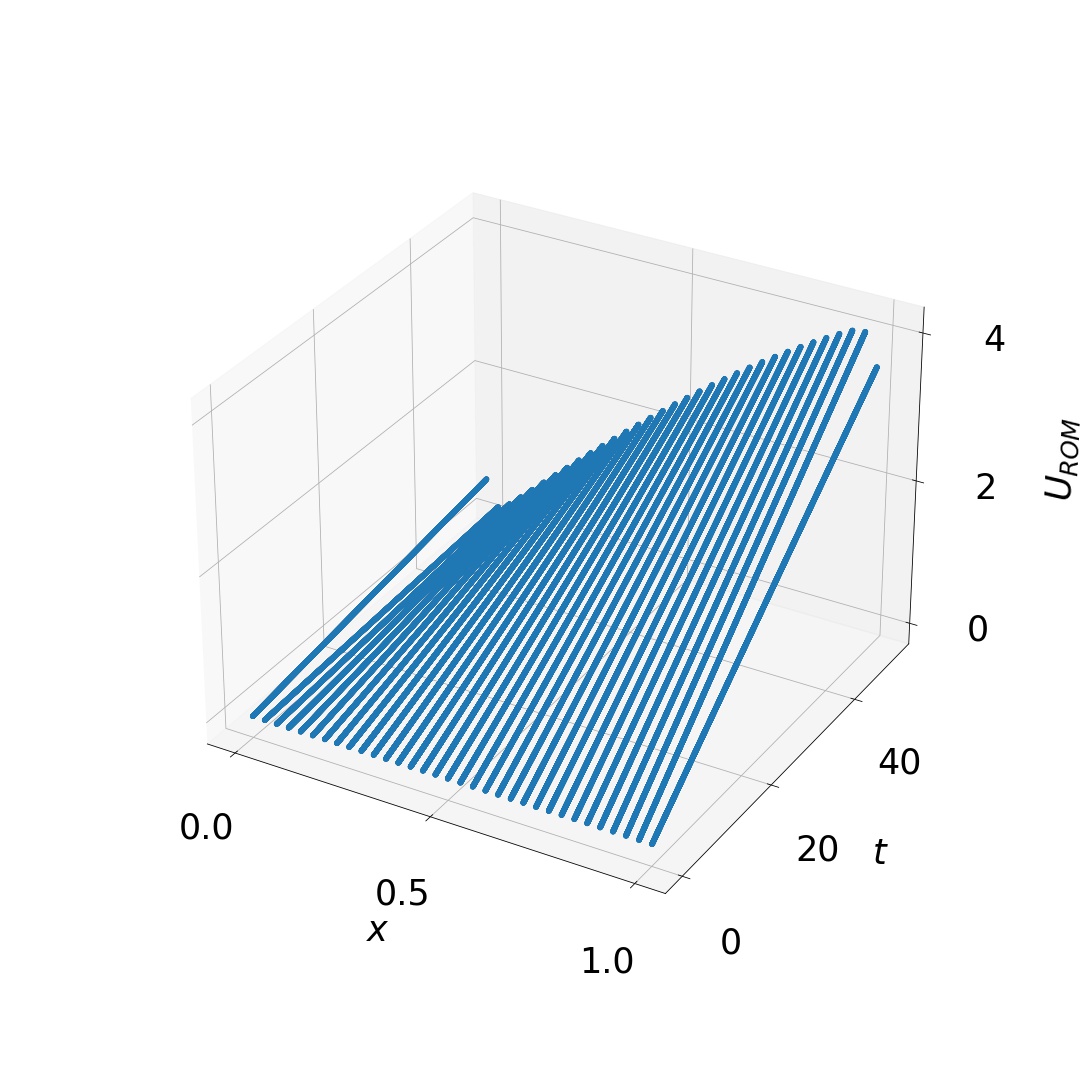}
		\caption{POD with $i=1$ modes.}
	\end{subfigure} \hspace{1em}
	\begin{subfigure}{0.4\linewidth}
		\includegraphics[width=\textwidth]{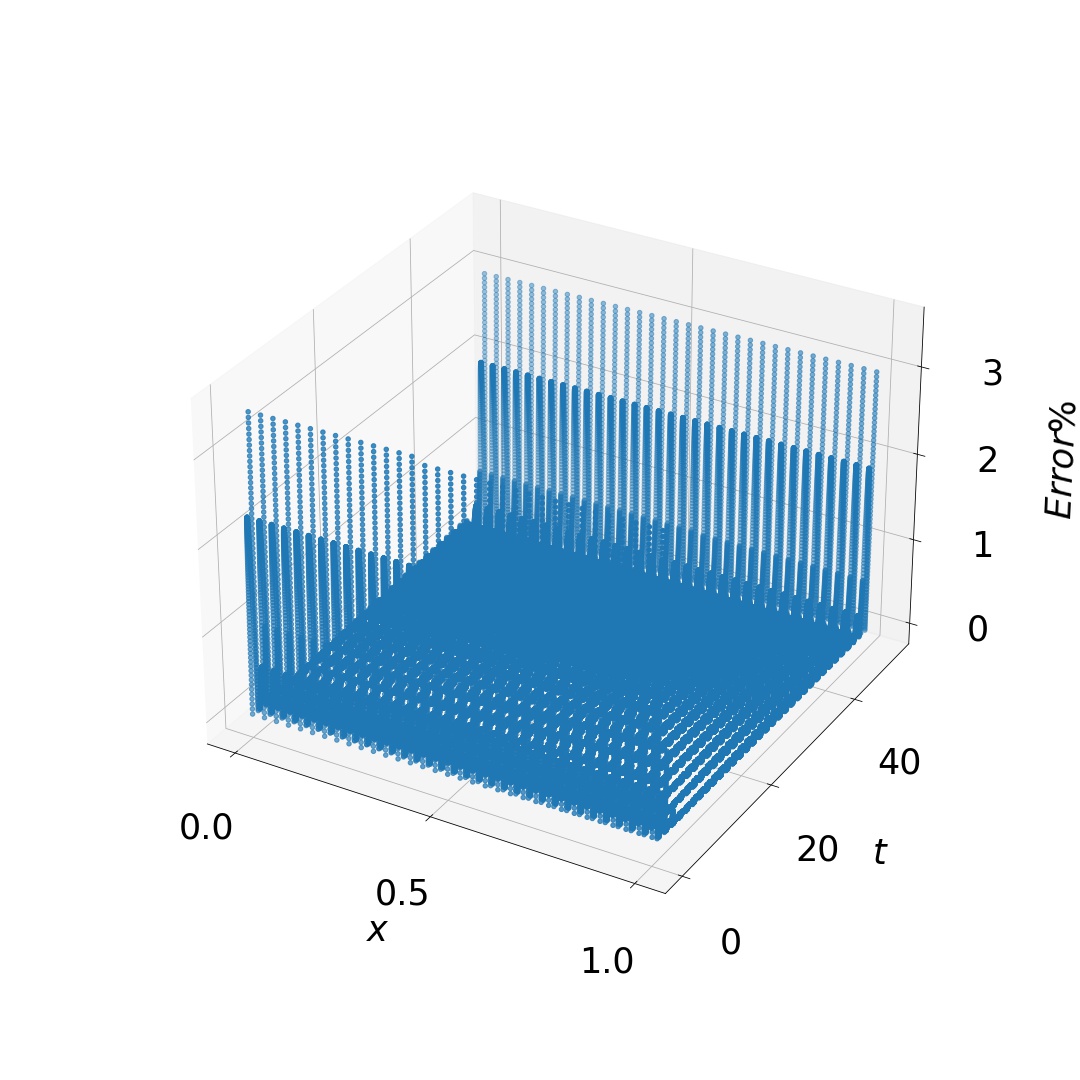}
		\caption{$\%$ Error with $i=1$ modes}
	\end{subfigure} \hspace{1em}
	\begin{subfigure}{0.4\linewidth}
		\includegraphics[width=\textwidth]{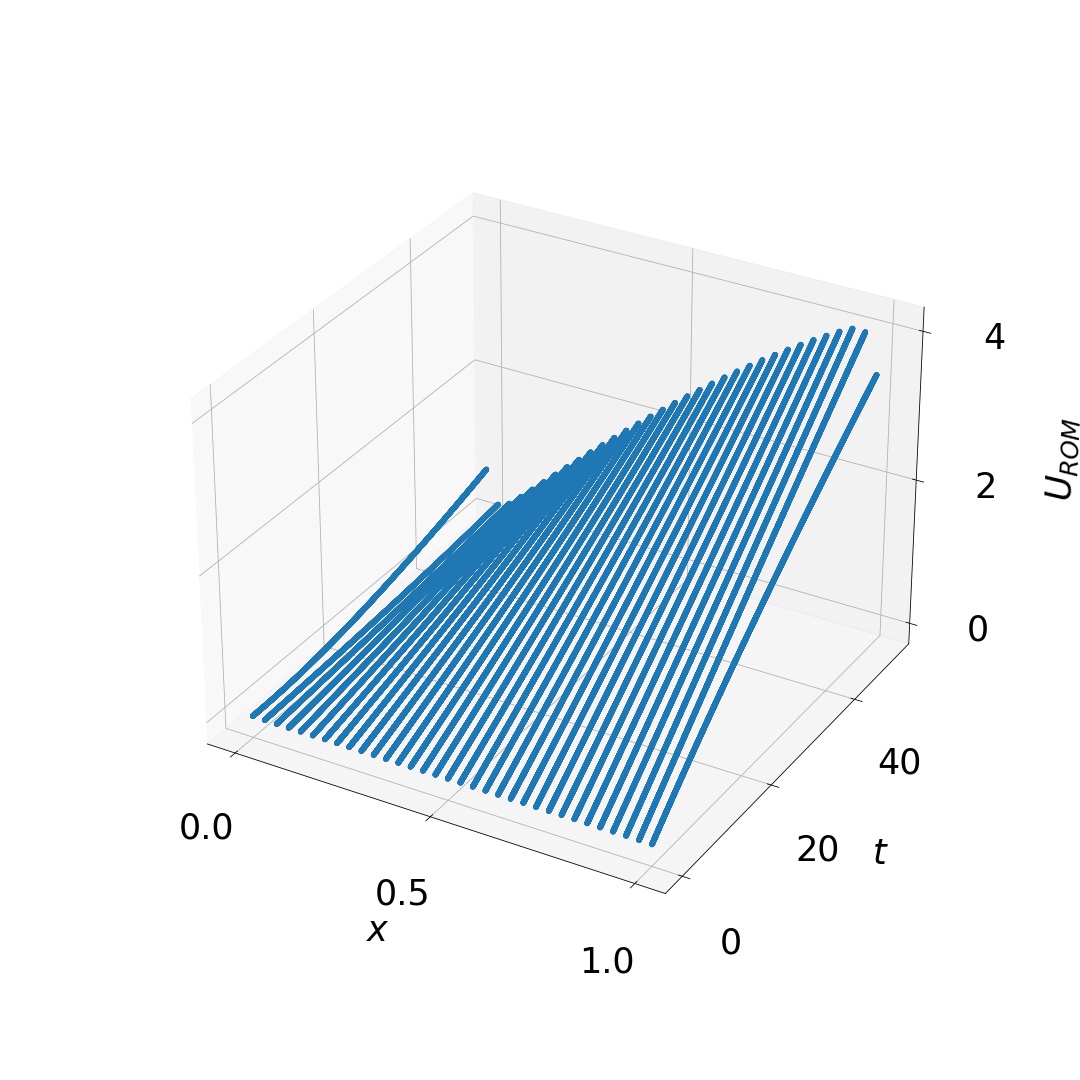}
		\caption{POD with $i=3$ modes.}
	\end{subfigure}
	\begin{subfigure}{0.4\linewidth}
		\includegraphics[width=\textwidth]{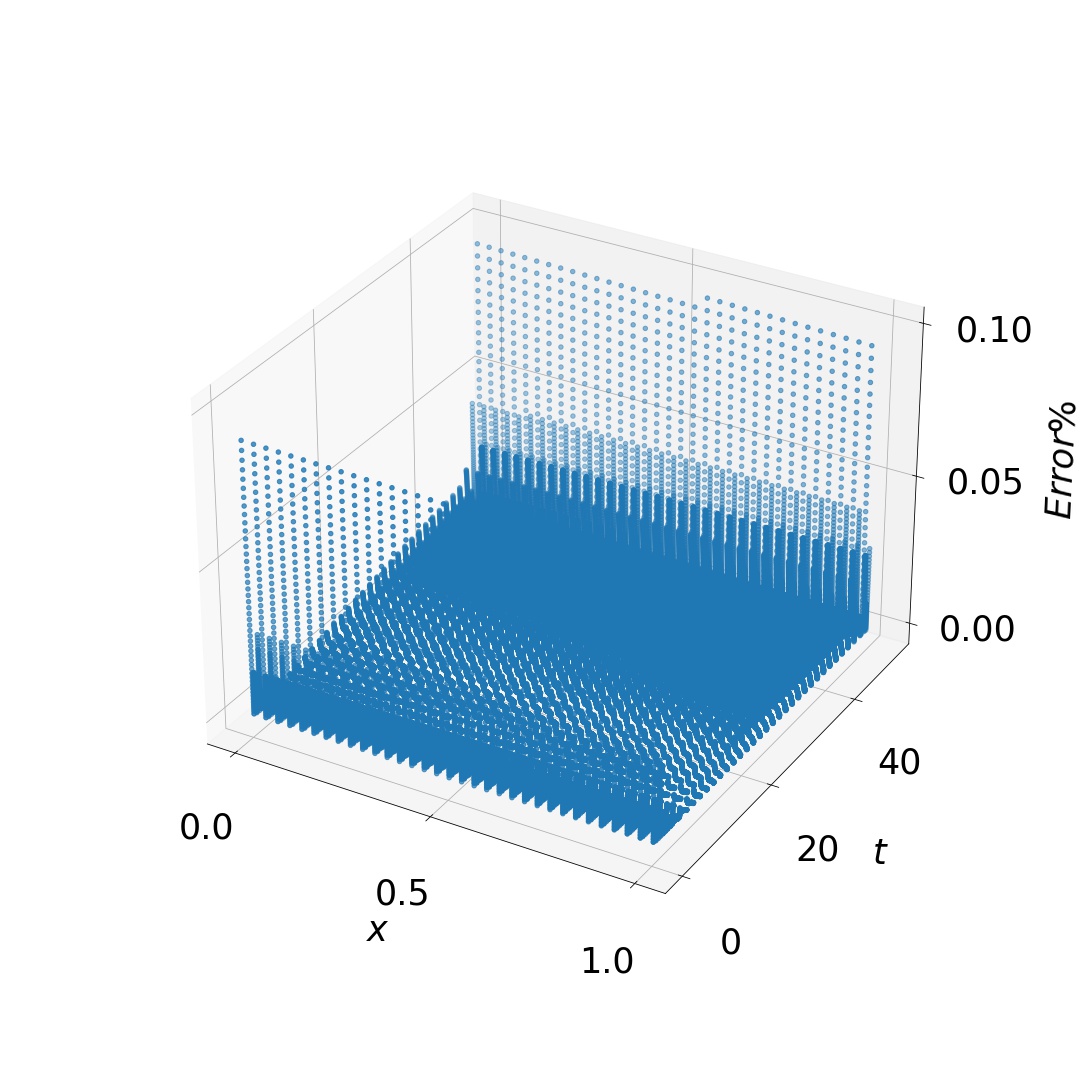}
		\caption{$\%$ Error with $i=3$ modes}
	\end{subfigure} \hspace{1em}
	\caption{Case 5 - Finite difference simulation of 1D Heat equation and vs POD-ROM}
	\label{fig:gpod_case5}
\end{figure}

Snapshots of the field variable, $U \in \mathbb{R}^{m \times n}$ of $m$ spatial points and $n$ temporal points are obtained from the numerical simulation of 1D Burgers' equation. Singular value decomposition of this matrix, decomposes the solution into,
\begin{equation}
    SVD(U) = \Phi \Sigma \Psi^*,
\end{equation}
where, the matrix $\Phi \in \mathbb{R}^{m \times m}$ represents the left singular vector and $\Psi^* \in \mathbb{R}^{n \times n}$ the conjugate transpose of the right singular vector. These matrices are composed of the eigenvectors $\phi$ (columns of $\Phi$), $\psi$ (rows of $\Psi^*$) of the covariance matrices $U^*U$ and $U U^*$. The matrix $\Sigma \in \mathbb{R}^{m \times 1}$ represents the singular values of $U$. The singular values are ordered as $\sigma_1 \ge \sigma_2 \ge...\ge \sigma_m$. The $i^{th}$ row of the solution matrix thus can be written as,
\begin{equation}
    u_{i} = \sum_{i} ^m \sigma_i \phi_{ij} \psi_{ik}.
\end{equation}
$\sigma_i$ represents the magnitude of significance of the modes. The relative significance of the modes can be estimated as,
\begin{equation}
    r = \frac{\sigma_i}{\sum_{i=1} \sigma_i},
\end{equation}
where, $\sum r=1$. As observed in Figs. \ref{fig_energy_case1}-\ref{fig_energy_case3}, when the diffusivity coefficient ($\nu$) is high, the first three to five modes contain upto $99 \%$ of the energy. A reduced solution can then be obtained by choosing a few dominant modes ($r$). Once the low dimensional basis set has been determined, it can be used to reformulate the Burgers' equation to create a reduced order model. Separation of variable and basis expansion is applied to Eq.\ref{eq_burgers}, where the field variable is expressed as
\begin{equation}
    U(x,t)= A(t) \Phi(x) = \sum_{i=1} ^m a_{ik} (t) \phi_{ij}(x),
    \label{sepvar}
\end{equation}
where, $\phi_{ij}(x)$ represents the orthonormal basis functions obtained by SVD of the data matrix, $a_{ik}(t)$ the temporal coefficients, and $i$ index of number of modes. It is to be noted that the basis function and temporal coefficients exist in a function space, where the continuous function $U(x,t)$ is discretized into a system of finite dimensions. By substituting Eq.\ref{sepvar} into Eq.\ref{eq_burgers}, the PDE for $U(x,t)$ can be transformed into an ordinary differential equation (ODE) for $A(t)$. 
\begin{equation} \label{eq_burger_sep}
\begin{aligned}
    \frac{d}{dt}\sum a_{ik} \phi_{ij} = - \sum a_{ik} \phi_{ij} \frac{d}{dx} \sum a_{ik} \phi_{ij} + \nu \frac{d ^2}{dx^2} \sum a_{ik} \phi_{ij} + Q, \\
\end{aligned}
\end{equation}
Since the basis functions are orthogonal they hold the property,
\begin{equation}
    \int \phi_{ij}(x) \phi_{il}(x)^* dx =
    \sum \phi_{ij}(x) \phi_{il} (x) = 
    \begin{cases}
      1 & j=l\\
      0 & j \neq l
    \end{cases}
\end{equation}
The final form of equation is obtained by multiplying both sides of Eq.\ref{eq_burger_sep} by $\phi_{ij}^T(x)$, 
\begin{equation}
\label{eq_gpod}
\begin{aligned}
    \frac{d}{dt}\sum a_{ik} = - \phi_{ij}^T \sum a_{ik} \phi_{ij} \frac{d}{dx} \sum a_{ik} \phi_{ij} + \phi_{ij}^T \nu \frac{d ^2}{dx^2} \sum a_{ik} \phi_{ij} + \phi_{ij}^T Q.    
    \end{aligned}
\end{equation}
To build a reduced order model, a small number of modes ($i<<m$) are considered in Eq.\ref{sepvar} and Eq.\ref{eq_gpod}. The spatial derivatives in Eq.\ref{eq_gpod} are solved using finite difference approach. Periodic boundary conditions are applied and the same initial condition (as the one used in the finite difference simulation), is imposed for the temporal coefficient $a_0=U(x,0) \Phi$. This equation is then solved for $a(t)$ using a standard ODE solver and a ROM is created by projection of this solution on $\Phi$. The results from the ROM are compared against the numerical simulation and presented in Figs. \ref{fig:gpod_case1}-\ref{fig:gpod_case5}. It is observed that the ROM is very accurate in predicting the field variable's evolution in space and time with a few modes ($i \le 5$). For new initial and boundary conditions, it is sufficient to assemble a snapshot matrix considering different instances of initial and boundary conditions and compute the POD at an offline stage. The ratio of computation time required to run the finite difference simulation and the POD-ROM for various cases are presented in Table \ref{tab:comptime}. In most cases, POD-ROM requires a smaller fraction of computation time when compared to the finite difference simulation. The costs of determination of POD basis and Galerkin projection increases with increase in number of dimensions, resolution and simulation run times.
\subsection{Implementation of POD}
\begin{outline}
    \1 Compute the instantaneous field variable $(U(x,t) \in \mathbb{R}^{m \times n})$ from the finite difference solution of 1D Burgers' equation.
    \1 Decompose the instantaneous field variable via Singular value decomposition,
    $$SVD(U)=U_L \cdot S \cdot V_\mathbb{R}^*.$$
    $U_L \in \mathbb{R}^{m \times m}$: Spatial modes (left singular vector), $S \in \mathbb{R}^{m \times 1}$: Magnitude of the modes (singular values), $V_\mathbb{R}^* \in \mathbb{R}^{n \times n}$: Temporal coefficients (complex conjugate of the right singular vector). Note: The singular values and vectors are ordered in a descending order.
    \1 Determine the number of significant modes by computing the relative magnitude of the singular values,
    $$r_k=\frac{S_k^2}{\sum_k S_k^2},k=1,2,\ldots,m.$$
    \1 Choose a small number of modes $(i)$ based on the values of $r (i \ll m).$
    \1 	Create new matrices with the reduced number of spatial modes and their corresponding temporal coefficients.
    $$\phi=U_L \in \mathbb{R}^{m \times i},$$
    $$\psi=V_\mathbb{R}^T \in \mathbb{R}^{i \times n},$$
    $$\sigma= S \in \mathbb{R}^{i \times 1}.$$
    \1 Specify initial condition (same initial condition as the numerical simulation). The initial temperature is projected onto the reduced basis,
    $$a_0=\phi \cdot U(x,0).$$
    \1 Compute first derivative $\phi_x$ (upwind scheme) and second derivatives $\phi_{xx}$ (central finite difference) of the spatial modes.
    \1 Assemble the parameters: modes, temporal coefficients, derivatives, initial conditions, boundary conditions, source, and diffusivity coefficients
    $$Parameters=[\phi^T,\phi,\phi_x,\phi_{xx},Q,\nu].$$
    \1 Perform Galerkin projection by applying separation of variable to the data matrix and combining it with the governing equation,
    $$U(x,t)=A(t) \Phi(x),$$
    $$\frac{d}{dt} \sum_i \phi_i a_i = -\sum_i \phi_i a_i \frac{d}{dx} \sum_i \phi_i a_i + \nu \frac{d^2}{dx^2} \sum_i \phi_i a_i + Q.$$
    \2 Invoke the orthogonality property,
     $$
     \int \phi_i \phi_j dx = 
    \begin{cases}
      1, & i=j\\
      0, & i \neq j
    \end{cases}
    $$
    \2 Multiply both sides of the equation by $\phi_i ^T$
     $$\frac{d}{dt} \sum_i a_i = -\phi_i ^T \sum_i \phi_i a_i \frac{d}{dx} \sum_i \phi_i a_i + \phi_i ^T \nu \frac{d^2}{dx^2} \sum_i \phi_i a_i + \phi_i ^T Q.$$
     \2 Compute the right hand side of the equation: 
     $$RHS = rhs_C + \nu \cdot rhs_D + rhs_S,$$
     where,
     \3 $U=\phi \cdot a$, which is equivalent to $\sum_i \phi_i a_i$,
     \3 $rhs_C = -\phi^T \cdot U \cdot \phi_x \cdot a$, which is equivalent to $-\phi_i ^T \sum_i \phi_i a_i \frac{d}{dx} \sum_i \phi_i a_i$,
     \3 $\nu \cdot rhs_D = \phi^T \cdot \nu \cdot \phi_{xx}$, which is equivalent to $\nu \frac{d^2}{dx^2} \sum_i \phi_i a_i,$
     \3 $rhs_S = \phi^T Q$, which is equivalent to $\phi_i ^T Q.$
    \1 Compute the temporal coefficient using a standard ODE solver,
    $$a(t) = Standard.ODE.Solver(RHS,a_0,t,Parameters).$$
    \1 Compute the time evolution of the reduced order model of the system,
    $$U_{ROM} (t) = a(t) \phi.$$
\end{outline}
\section{Results}
\begin{figure}[p]
	\centering
	\begin{subfigure}{0.4\linewidth}
		\includegraphics[width=\textwidth]{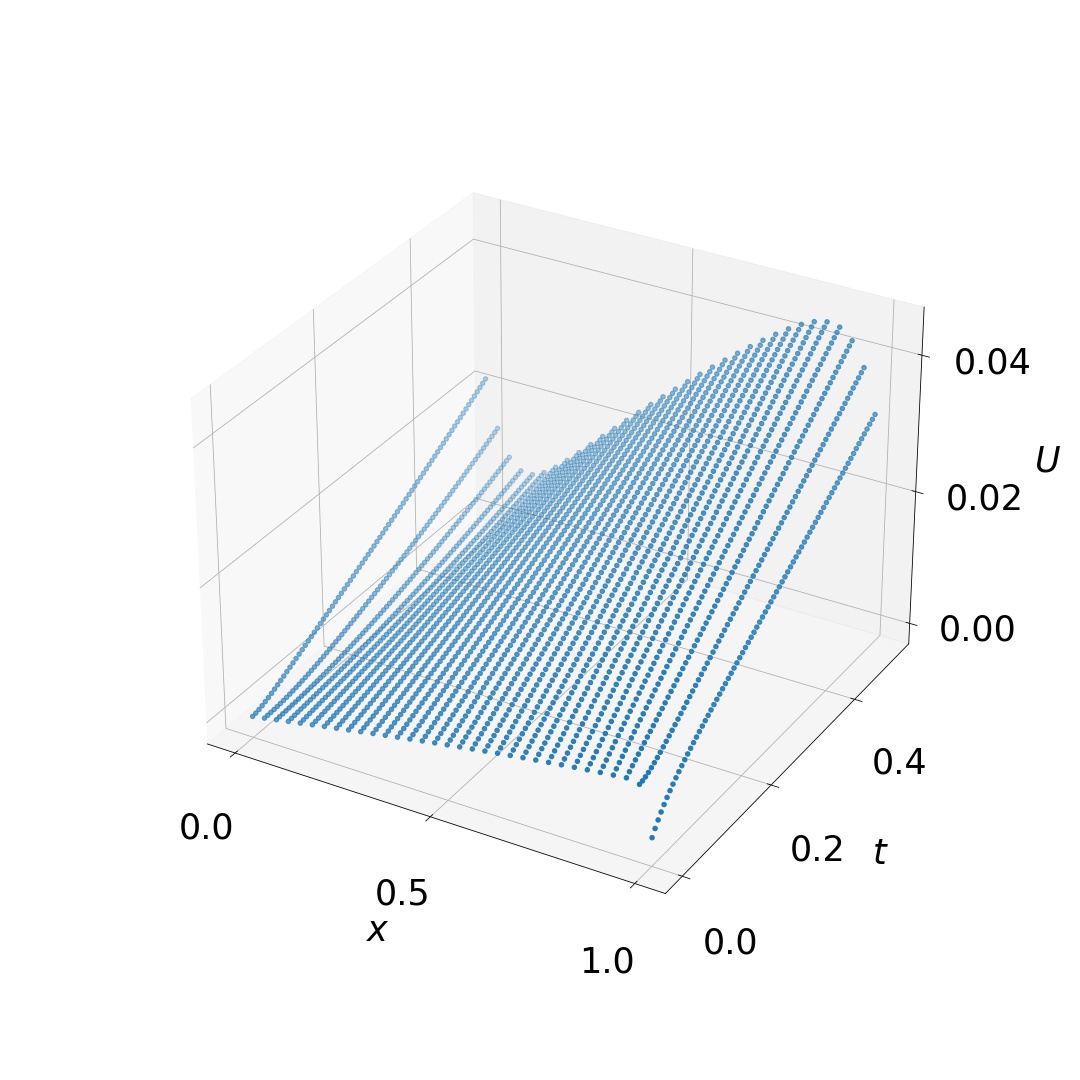} 
		\caption{Simulation with $t_{final}=0.5$, $N_x=32$, $CFL=0.2$, $\nu=0.01$. Run time $t_{Sim}=0.0615 s$}
		\label{fig_sim_case1}
	\end{subfigure} \hspace{1em}
	\begin{subfigure}{0.4\linewidth}
		\includegraphics[width=\textwidth]{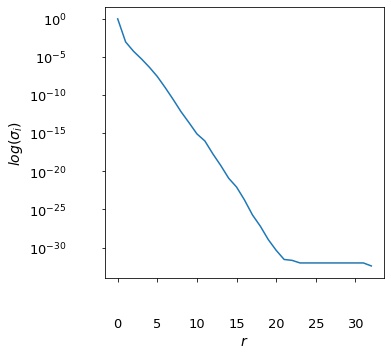} 
		\caption{Energy of dominant modes in the simulation}
		\label{fig_energy_case1}
	\end{subfigure} \\
	\begin{subfigure}{0.4\linewidth}
		\includegraphics[width=\textwidth]{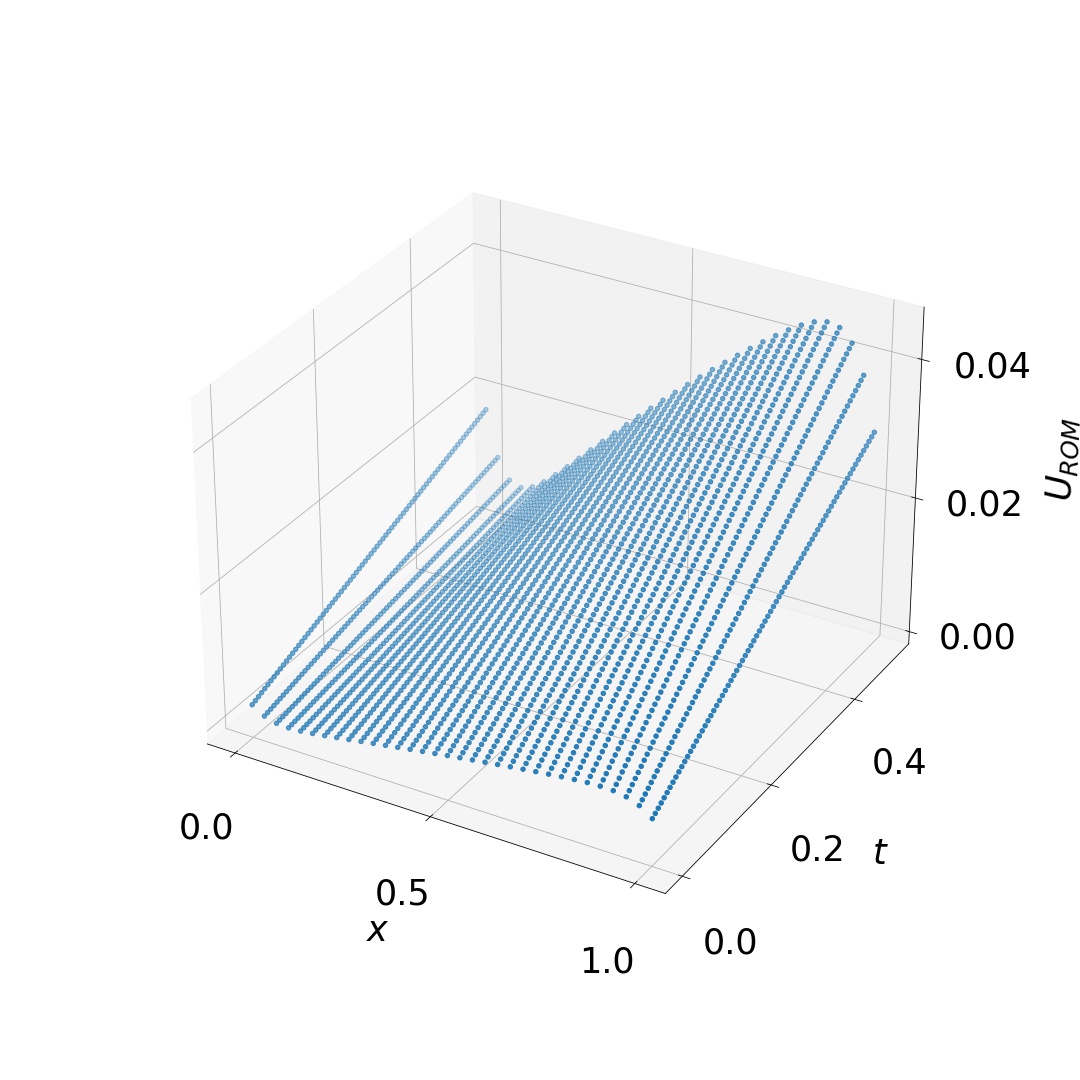}
		\caption{POD with $i=1$ mode. Model run time $t_{ROM}=0.0475 s$}
	\end{subfigure} \hspace{1em}
	\begin{subfigure}{0.4\linewidth}
		\includegraphics[width=\textwidth]{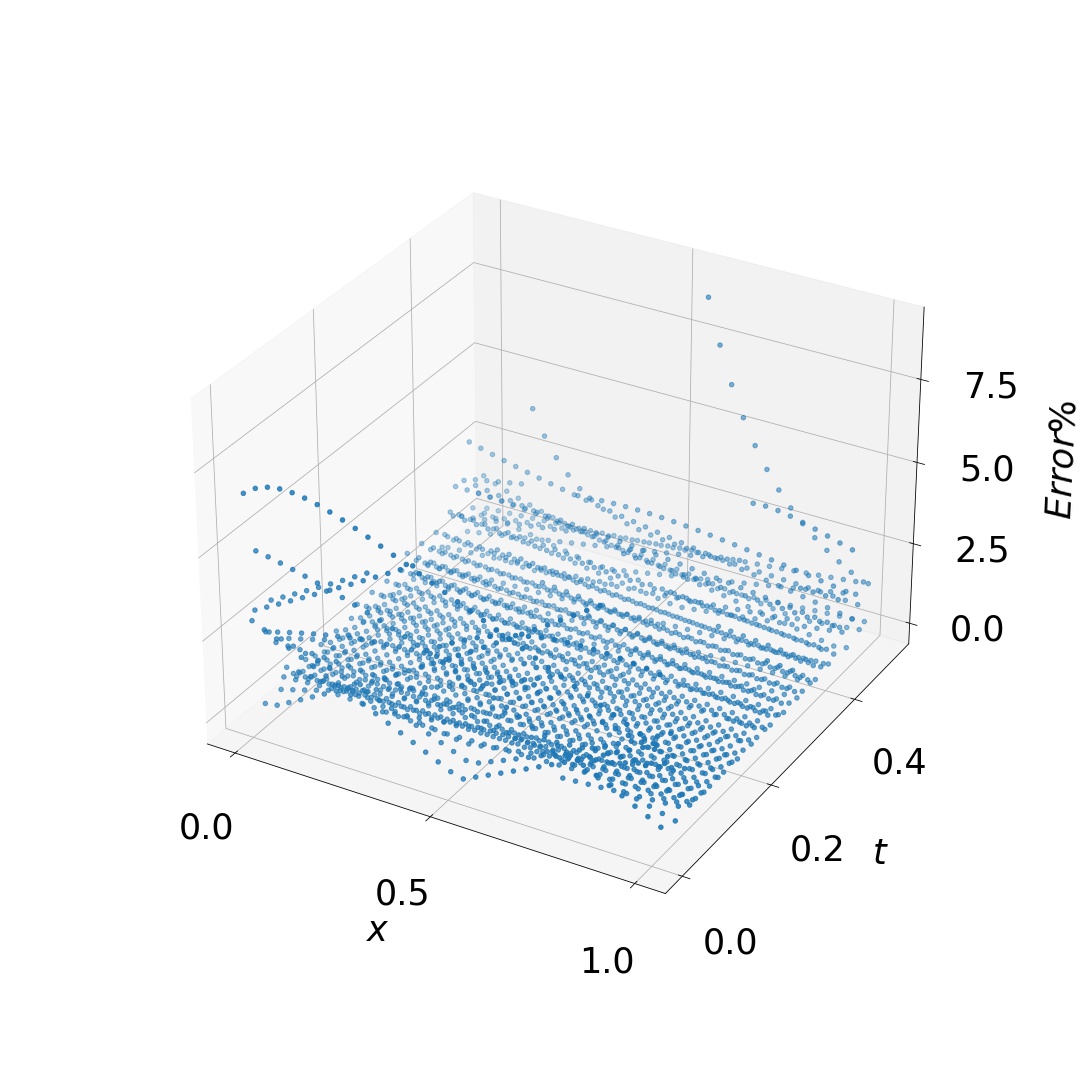}
		\caption{$\%$ Error with $i=1$ mode}
	\end{subfigure} \hspace{1em}
	\begin{subfigure}{0.4\linewidth}
		\includegraphics[width=\textwidth]{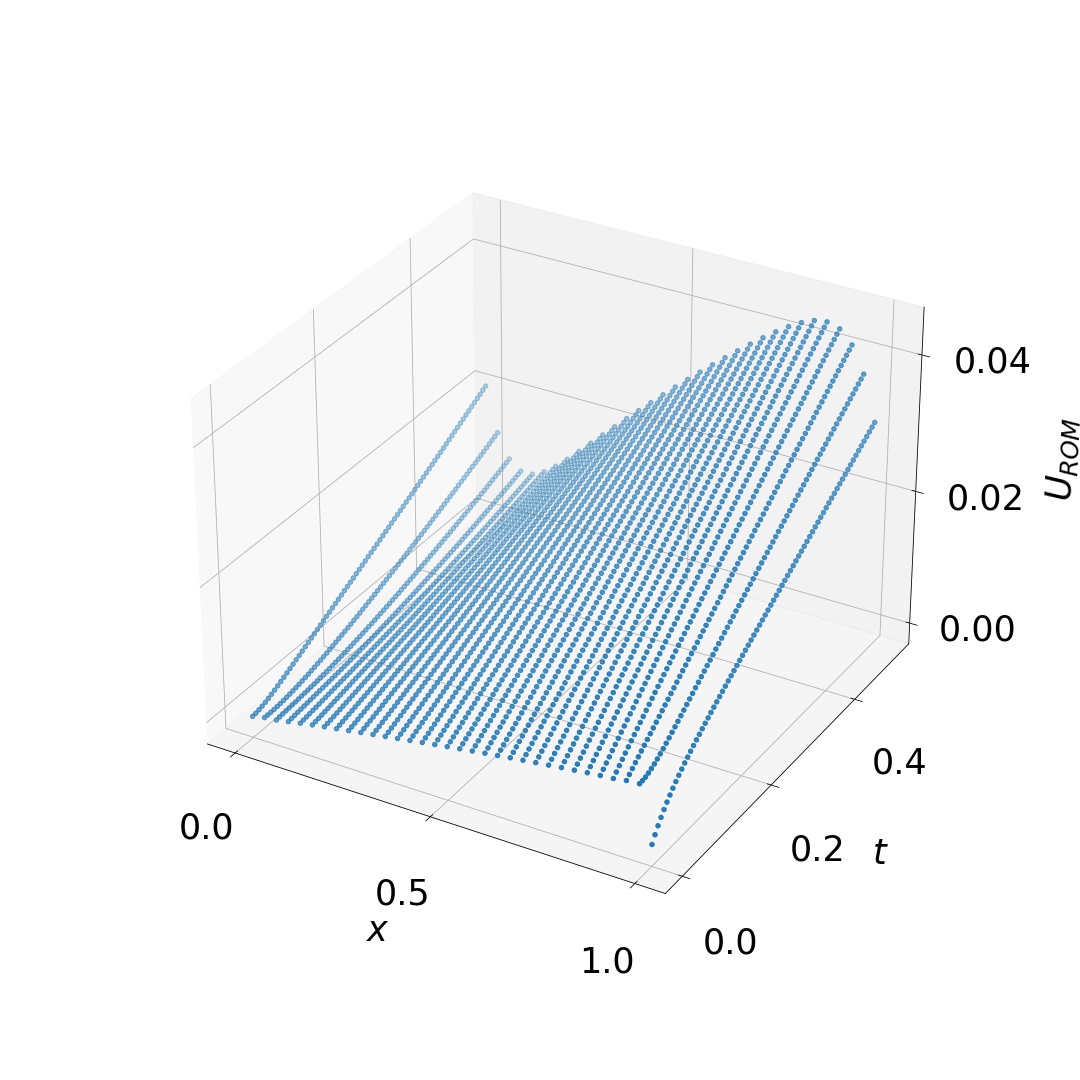}
		\caption{POD with $i=5$ modes. Model run time $t_{ROM}=0.064 s$}
	\end{subfigure}
	\begin{subfigure}{0.4\linewidth}
		\includegraphics[width=\textwidth]{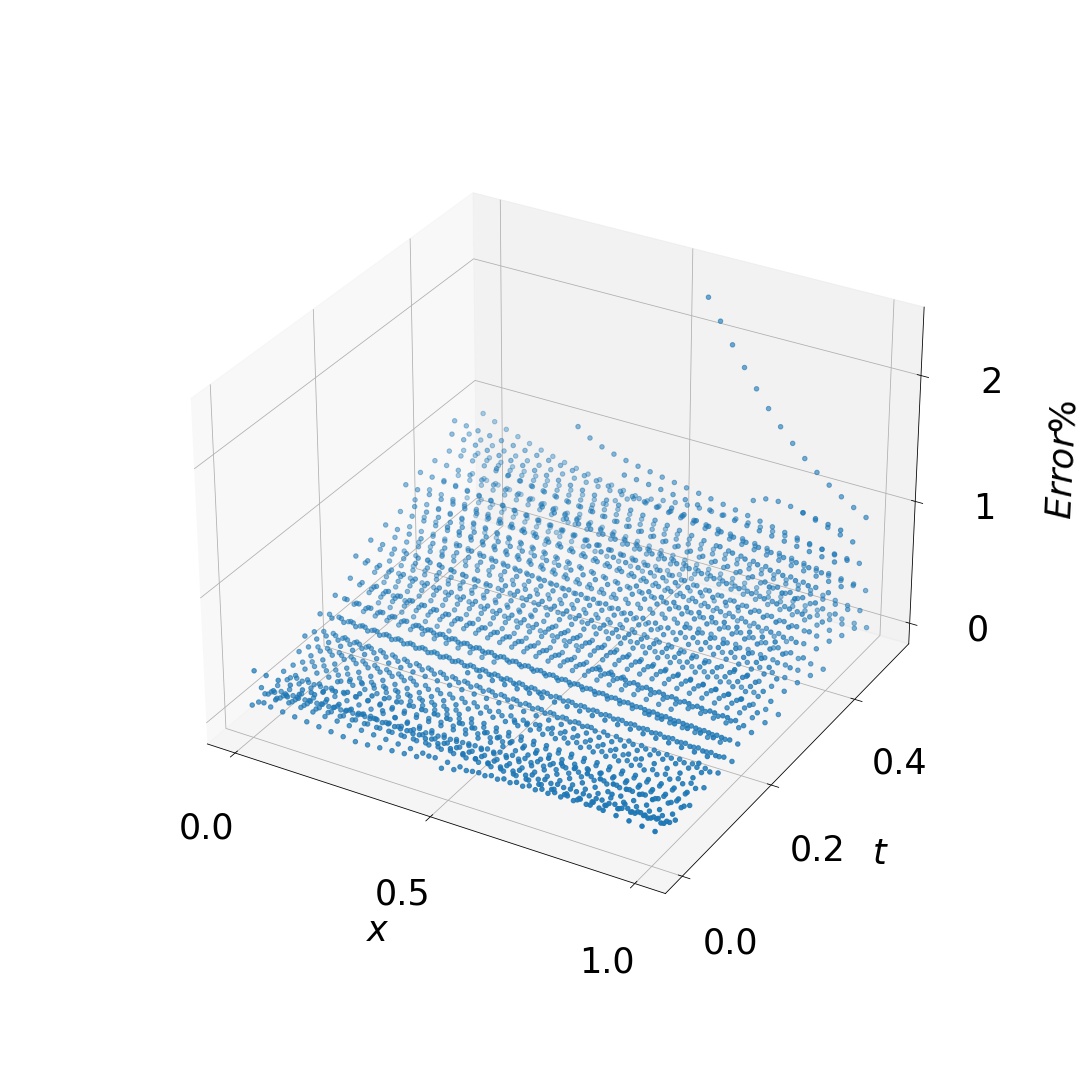}
		\caption{$\%$ Error with $i=5$ modes}
	\end{subfigure} \hspace{1em} 
	\caption{Case 1 - Finite difference simulation with short simulation time and high diffusivity vs reduced order model via Galerkin POD.}

	\label{fig:gpod_case1}
\end{figure}

\begin{figure}[!htb]
	\centering
	\begin{subfigure}{0.4\linewidth}
		\includegraphics[width=\textwidth]{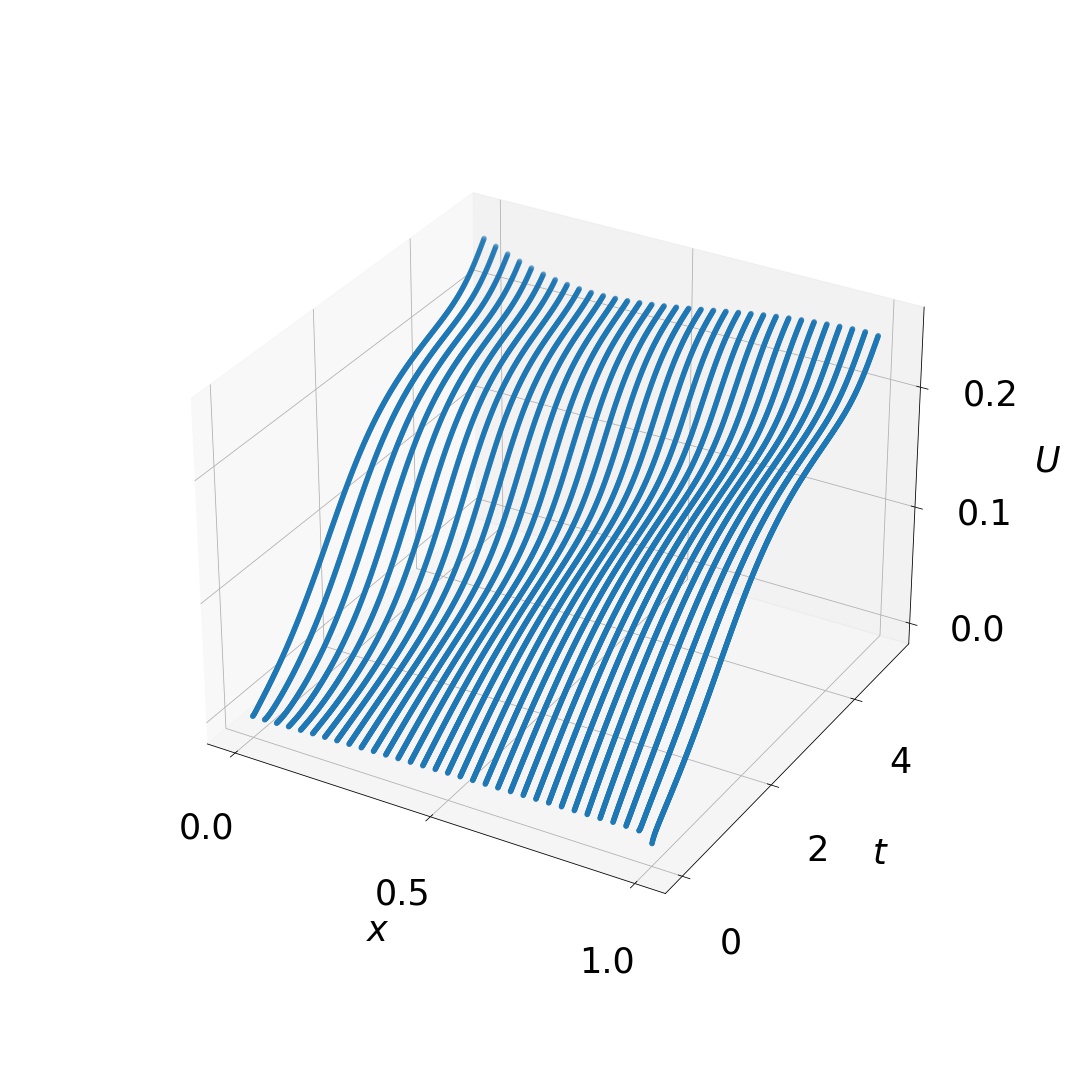} 
		\caption{Simulation with $t_{final}=5$, $N_x=32$, $CFL=0.2$, $\nu=0.01$. Run time $t_{Sim}=0.2031 s$}
		\label{fig_sim_case2}		
	\end{subfigure} \hspace{1em}
	\begin{subfigure}{0.4\linewidth}
		\includegraphics[width=\textwidth]{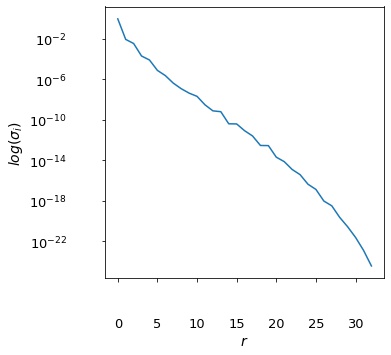} 
		\caption{Energy of dominant modes in the simulation}
		\label{fig_energy_case2}
	\end{subfigure} \\
	\begin{subfigure}{0.4\linewidth}
		\includegraphics[width=\textwidth]{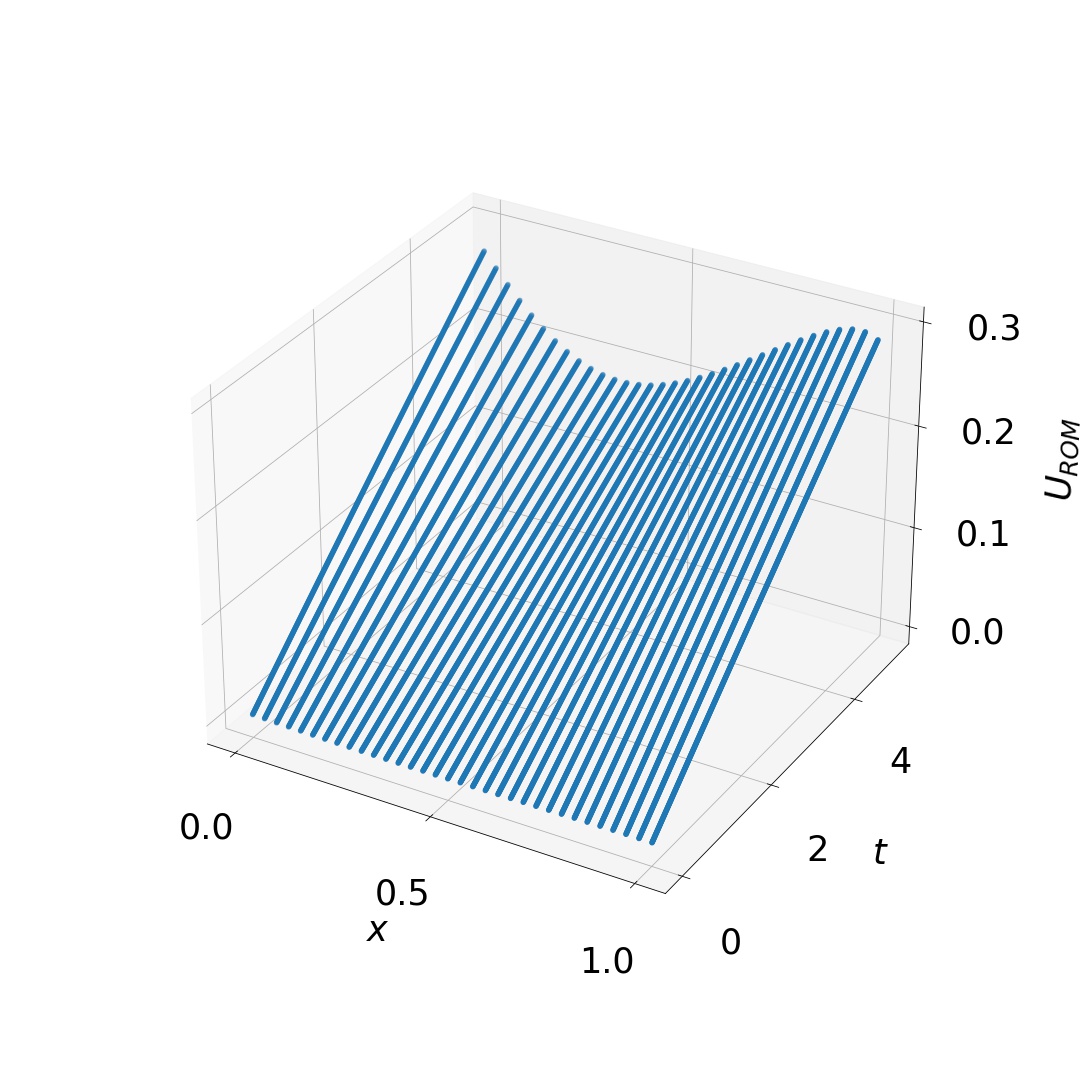}
		\caption{POD with $i=1$ mode. Model run time $t_{ROM}=0.0156 s$}
	\end{subfigure} \hspace{1em}
	\begin{subfigure}{0.4\linewidth}
		\includegraphics[width=\textwidth]{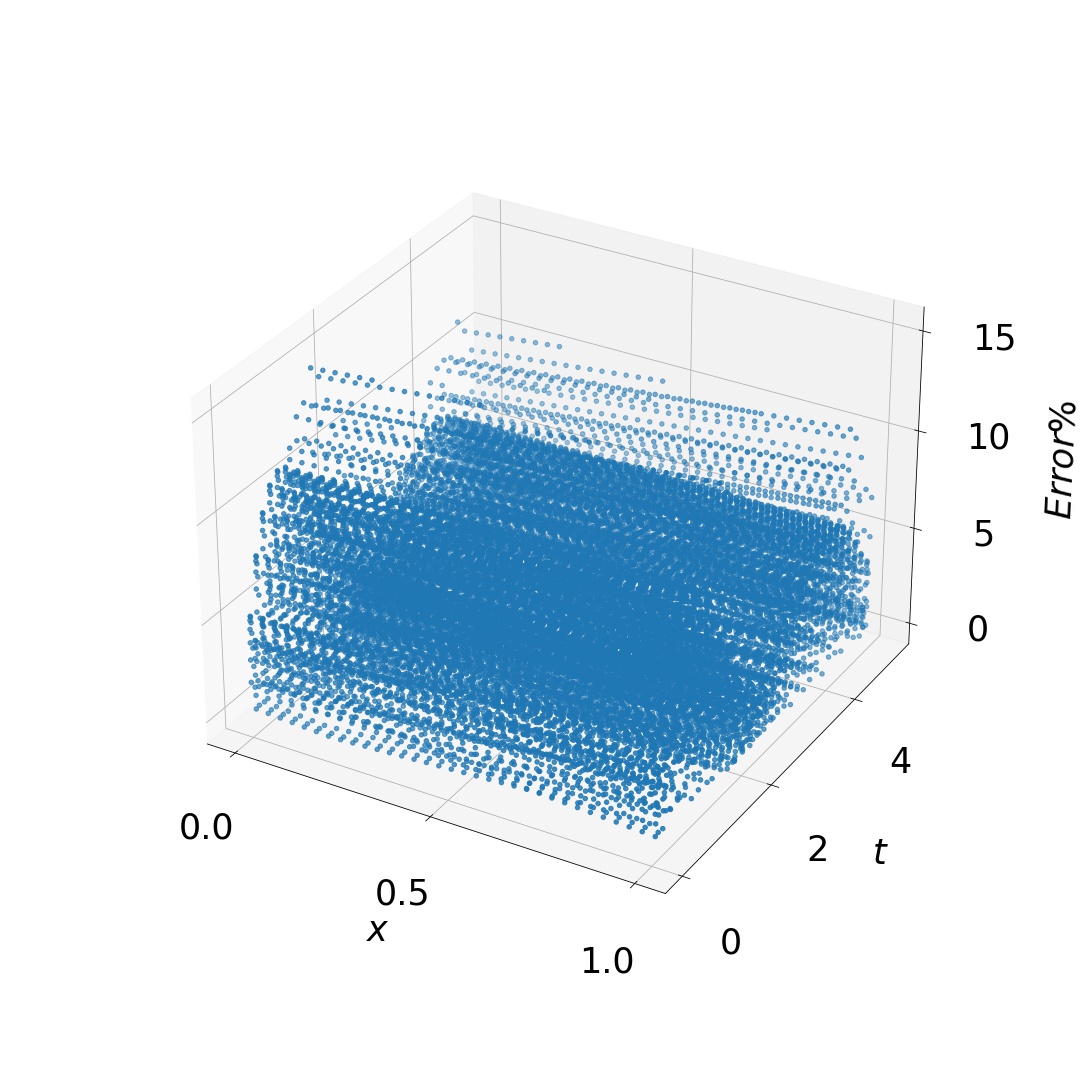}
		\caption{$\%$ Error with $i=1$ mode}
	\end{subfigure} \hspace{1em}
	\begin{subfigure}{0.4\linewidth}
		\includegraphics[width=\textwidth]{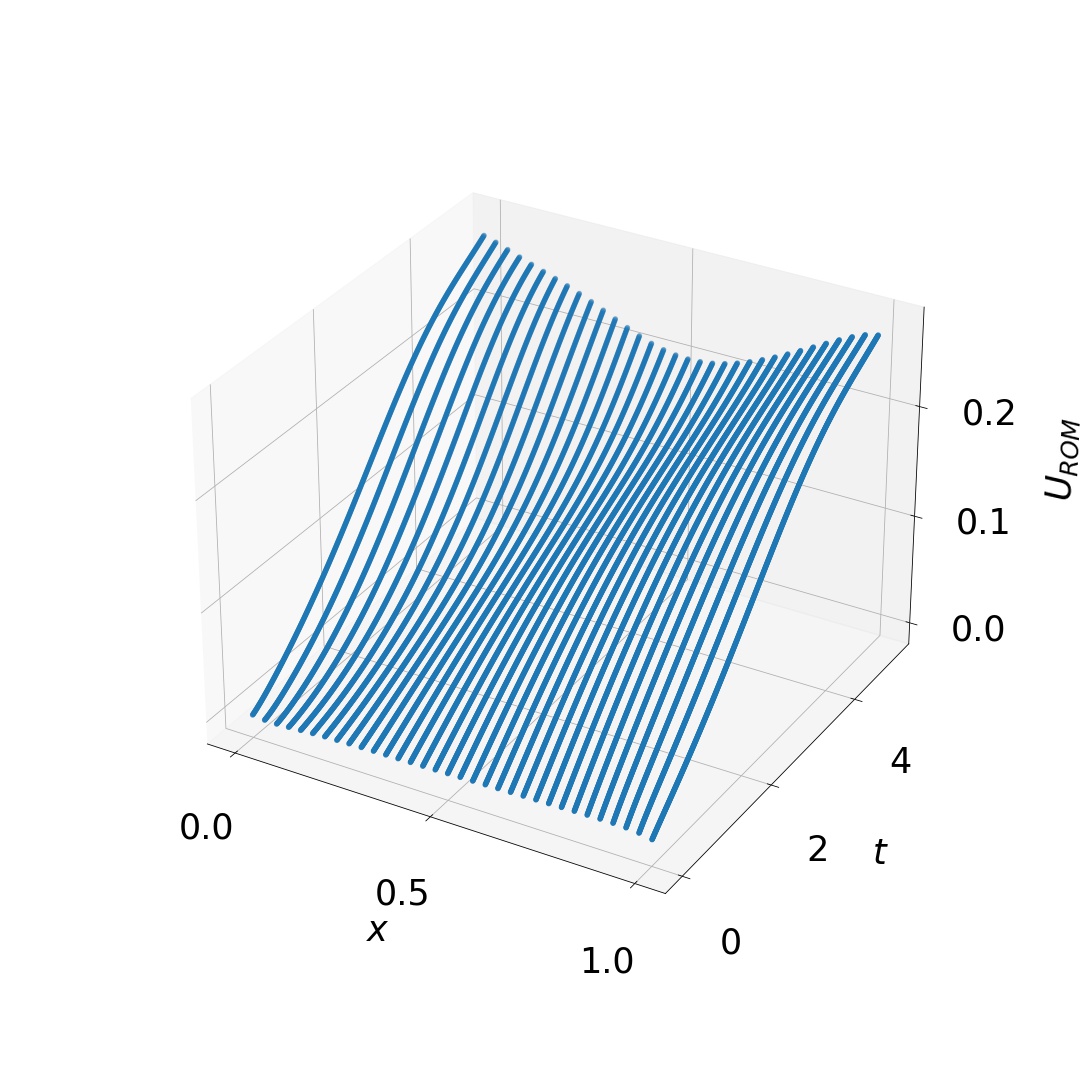}
		\caption{POD with $i=5$ modes. Model run time $t_{ROM}=0.031 s$}
	\end{subfigure}
	\begin{subfigure}{0.4\linewidth}
		\includegraphics[width=\textwidth]{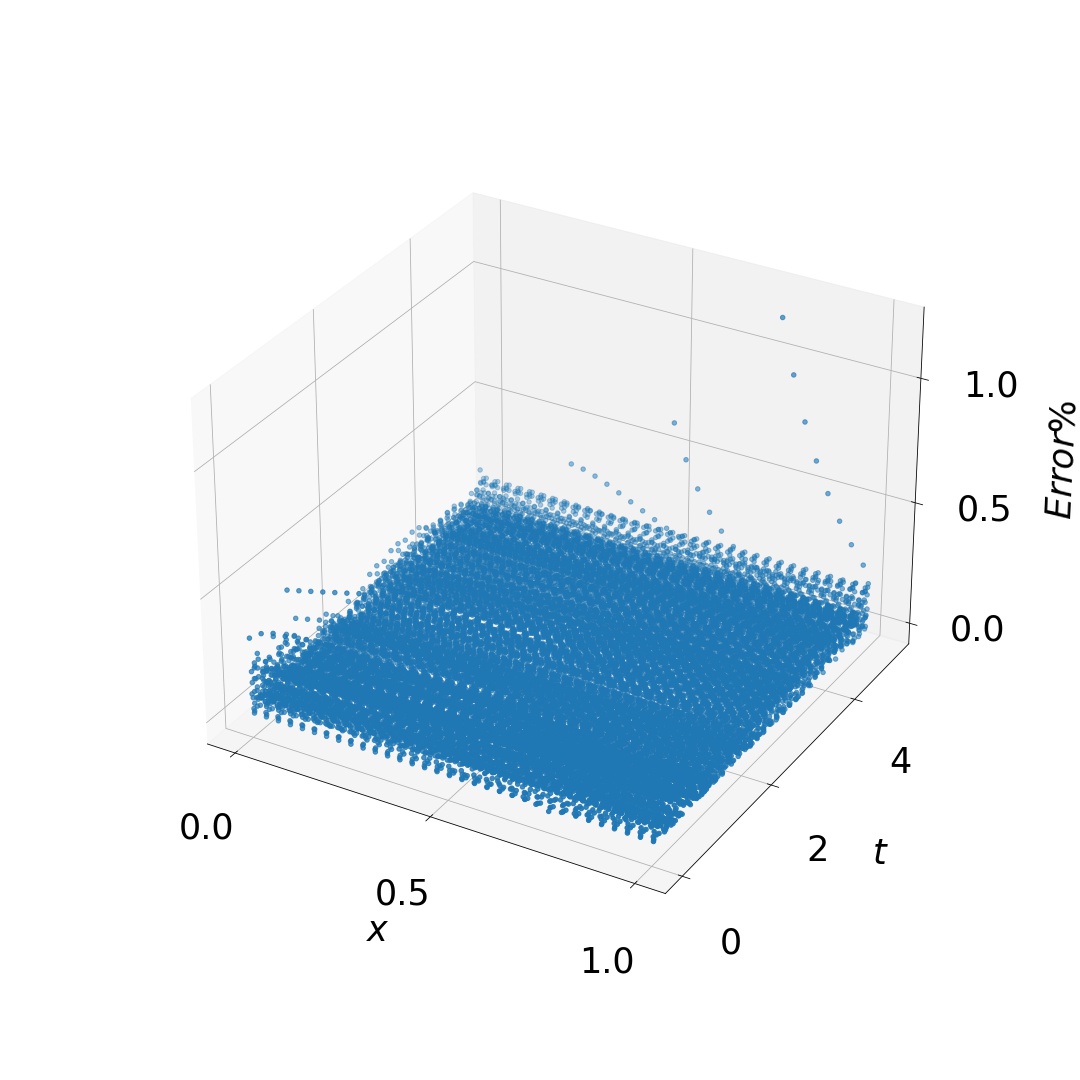}
		\caption{$\%$ Error with $i=5$ modes}
	\end{subfigure} \hspace{1em}
	\caption{Case 2 - Finite difference simulation with longer simulation time and high diffusivity vs reduced order model via Galerkin POD.}
	\label{fig:gpod_case2}
\end{figure}

\begin{figure}[!htb]
	\centering
	\begin{subfigure}{0.4\linewidth}
		\includegraphics[width=\textwidth]{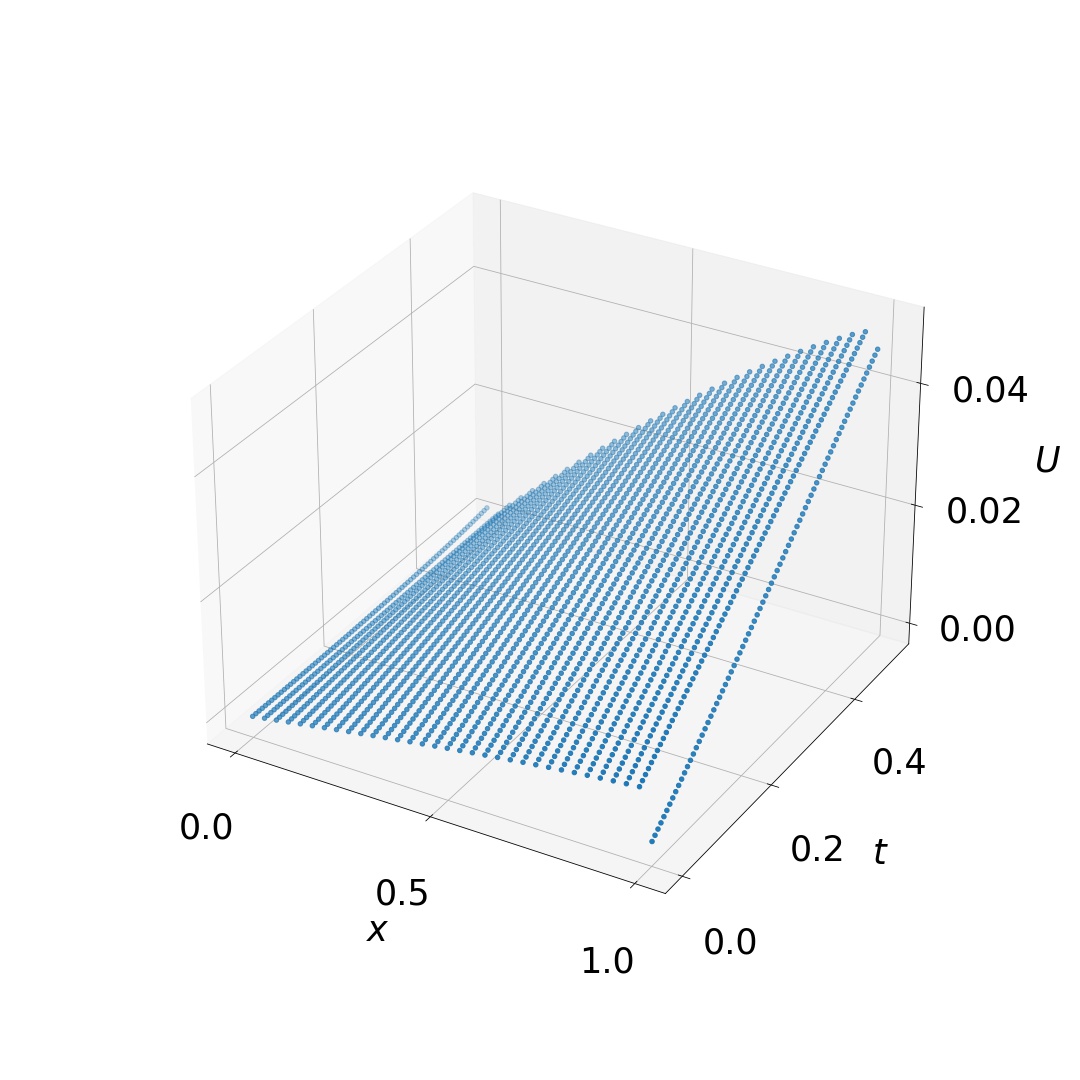} 
		\caption{Simulation with $t_{final}=0.5$, $N_x=32$, $CFL=0.2$, $\nu=0.0001$. Run time $t_{Sim}=0.049 s$}
		\label{fig_sim_case3}
	\end{subfigure} \hspace{1em}
	\begin{subfigure}{0.4\linewidth}
		\includegraphics[width=\textwidth]{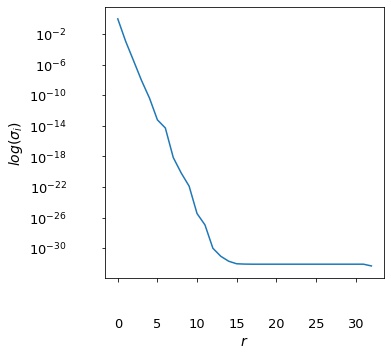} 
		\caption{Energy of dominant modes in the simulation}
		\label{fig_energy_case3}
	\end{subfigure} \\
	\begin{subfigure}{0.4\linewidth}
		\includegraphics[width=\textwidth]{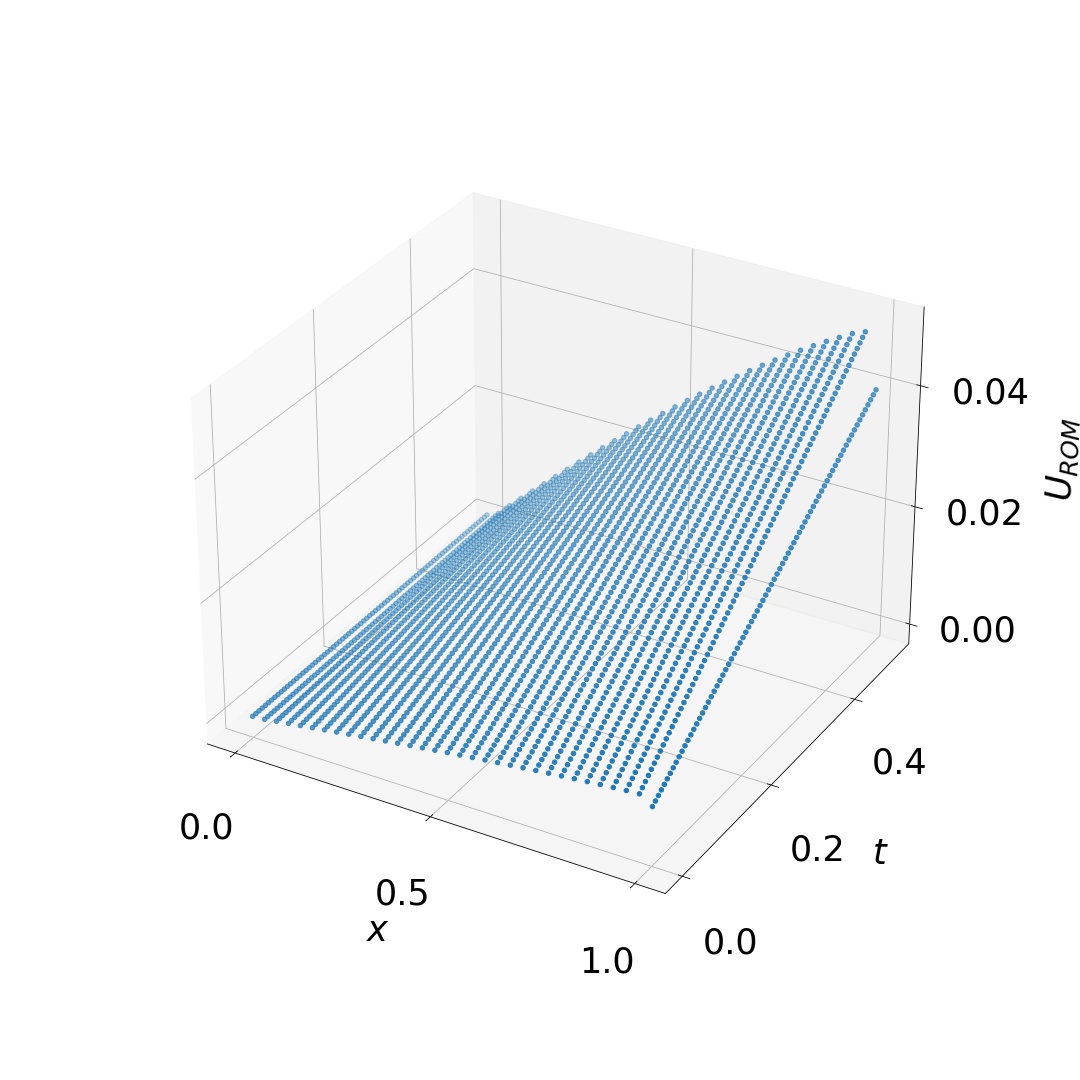}
		\caption{POD with $i=1$ mode. Model run time $t_{ROM}=0.028 s$}
	\end{subfigure} \hspace{1em}
	\begin{subfigure}{0.4\linewidth}
		\includegraphics[width=\textwidth]{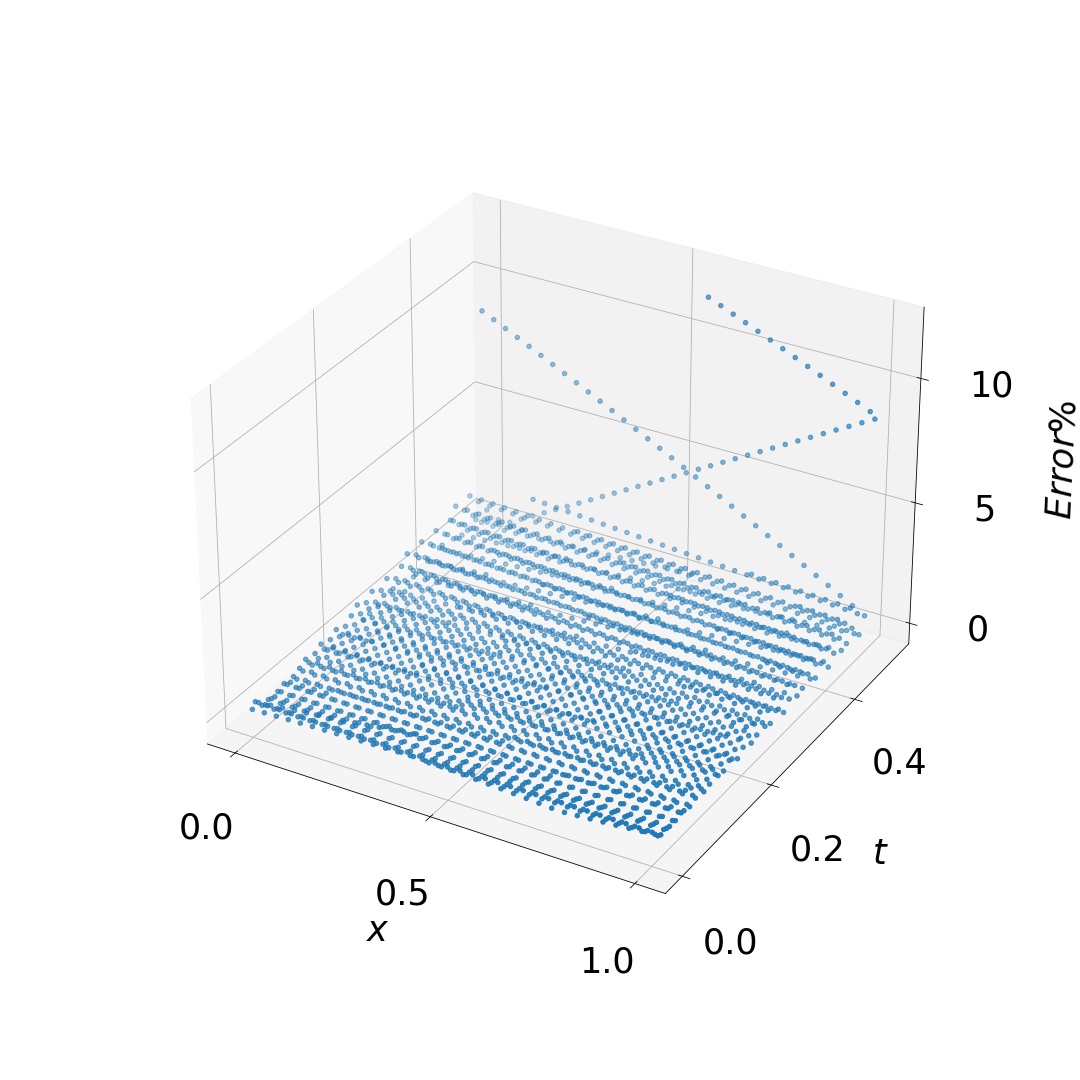}
		\caption{$\%$ Error with $i=1$ mode}
	\end{subfigure} \hspace{1em}
	\begin{subfigure}{0.4\linewidth}
		\includegraphics[width=\textwidth]{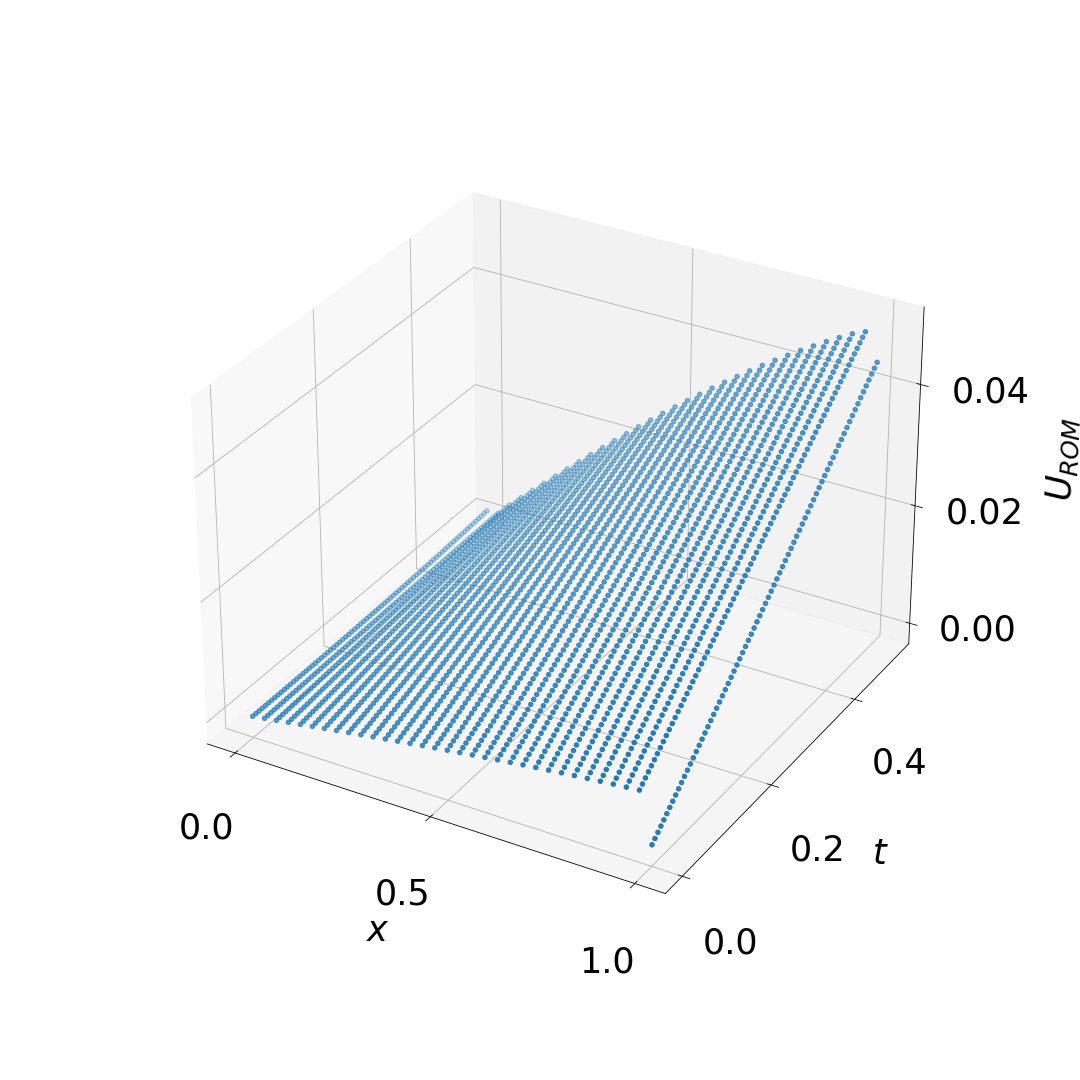}
		\caption{POD with $i=5$ modes. Model run time $t_{ROM}=0.054 s$}
	\end{subfigure}
	\begin{subfigure}{0.4\linewidth}
		\includegraphics[width=\textwidth]{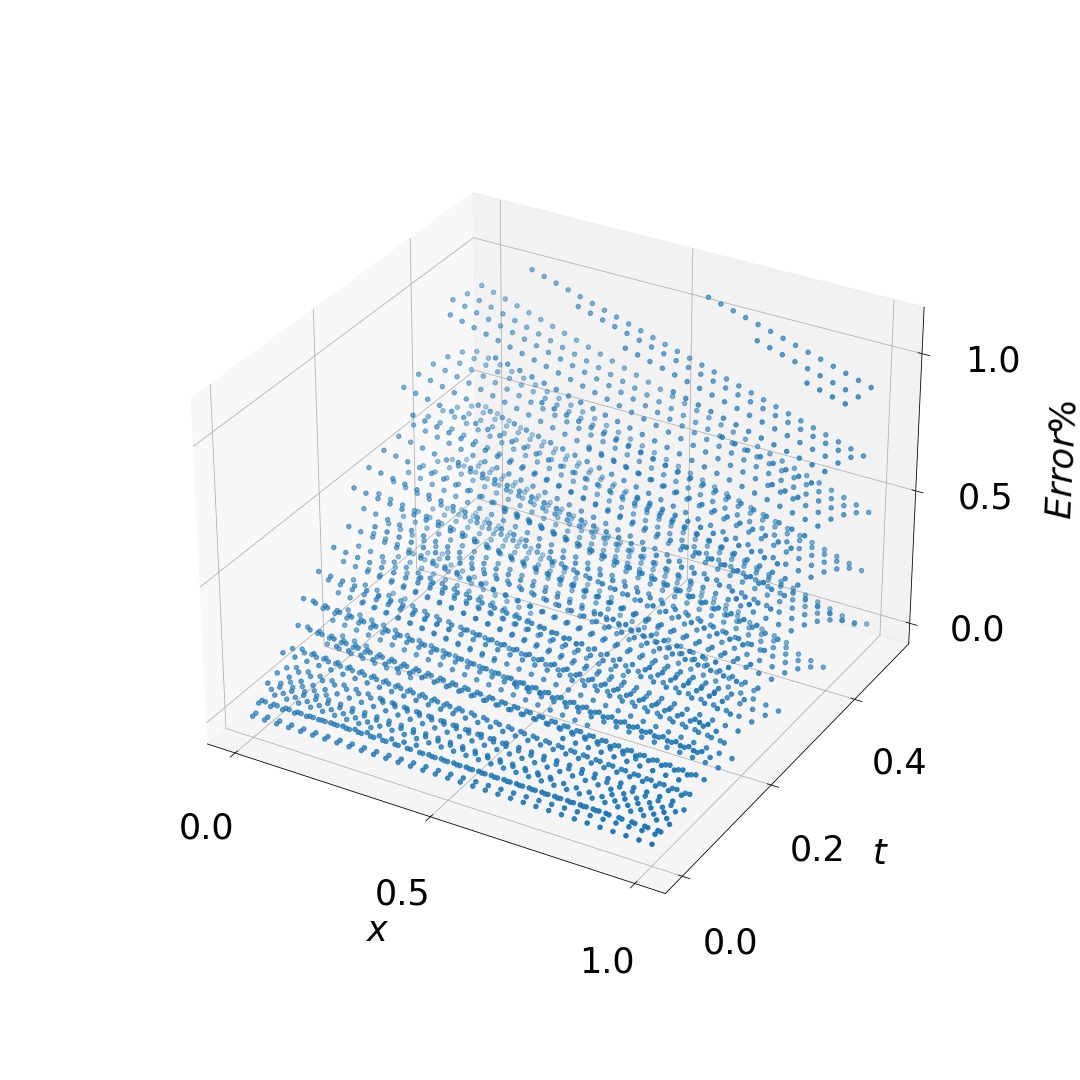}
		\caption{$\%$ Error with $i=5$ modes}
	\end{subfigure} \hspace{1em}
	\caption{Case 3 - Finite difference simulation with low diffusivity and short simulation time vs reduced order model via Galerkin POD.}
	\label{fig:gpod_case3}
\end{figure}

\begin{figure}[!htb]
	\centering
	\begin{subfigure}{0.4\linewidth}
		\includegraphics[width=\textwidth]{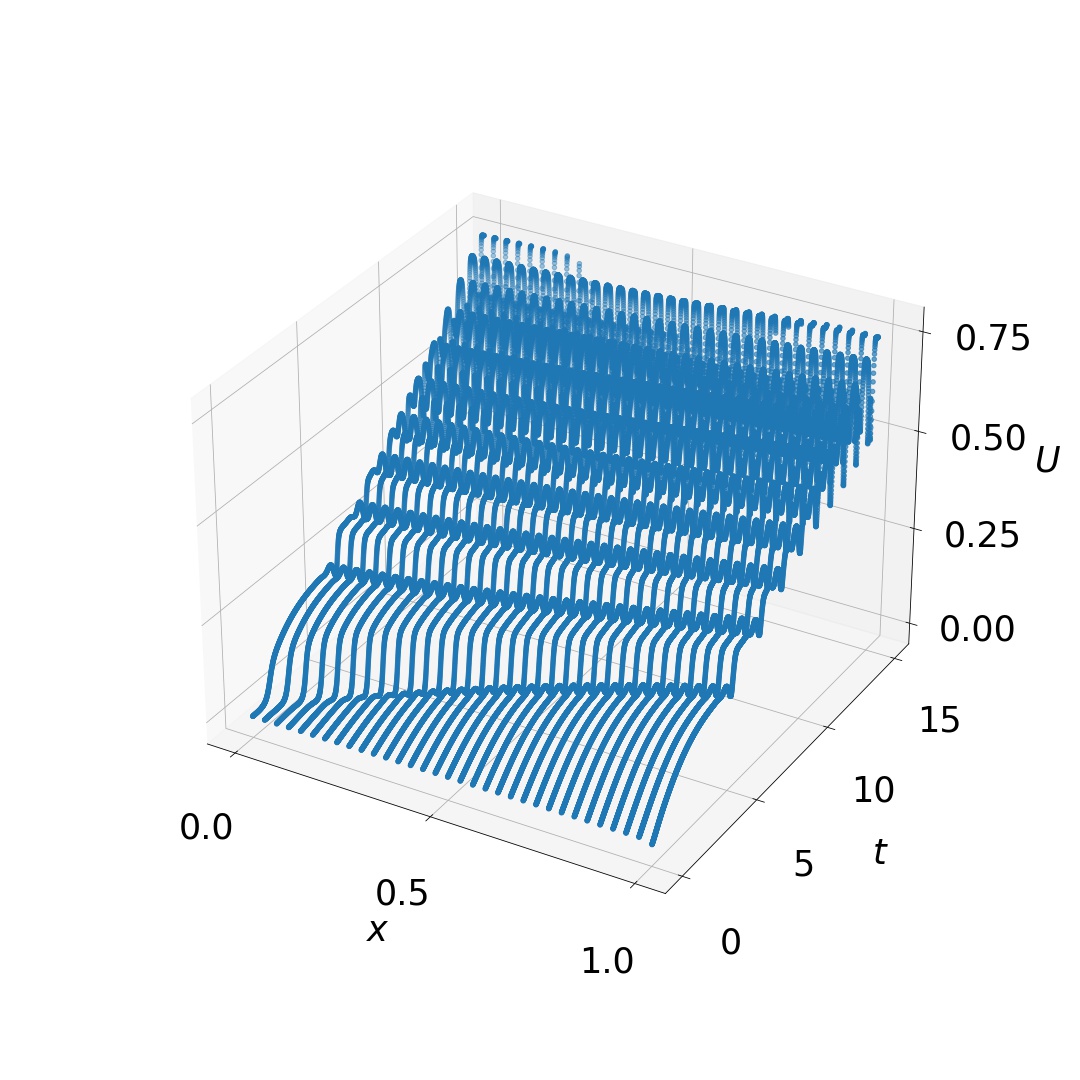} 
		\caption{Simulation with $t_{final}=15$, $N_x=32$, $CFL=0.2$, $\nu=0.0001$. Run time $t_{Sim}=0.648 s$}
		\label{fig_sim_case4}
	\end{subfigure} \hspace{1em}
	\begin{subfigure}{0.4\linewidth}
		\includegraphics[width=\textwidth]{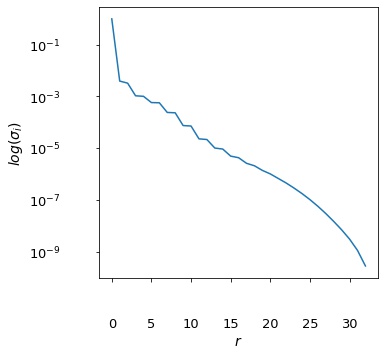} 
		\caption{Energy of dominant modes in the simulation}
		\label{fig_energy_case4}
	\end{subfigure} \\
	\begin{subfigure}{0.4\linewidth}
		\includegraphics[width=\textwidth]{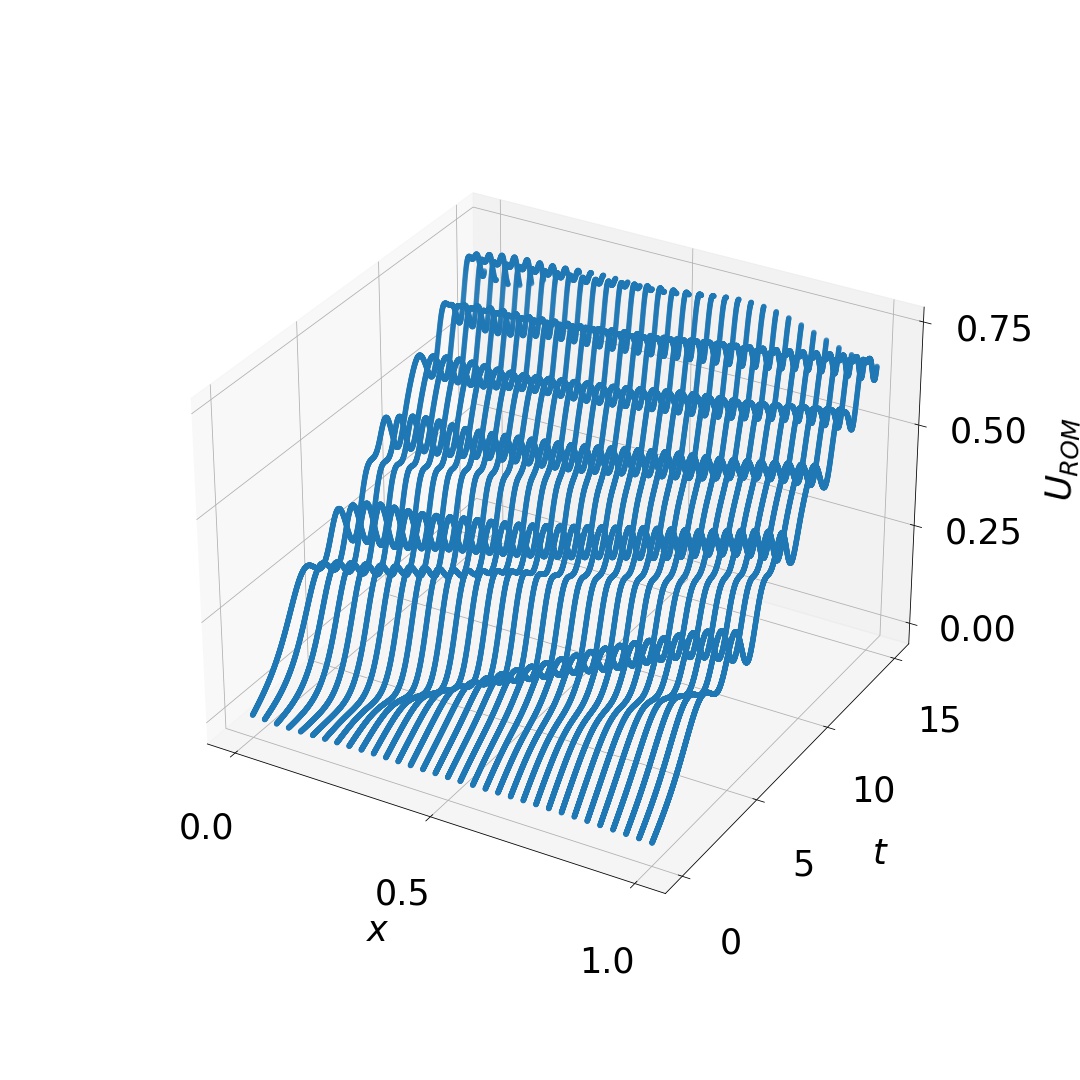}
		\caption{POD with $i=5$ modes. Model run time $t_{ROM}=0.203 s$}
	\end{subfigure} \hspace{1em}
	\begin{subfigure}{0.4\linewidth}
		\includegraphics[width=\textwidth]{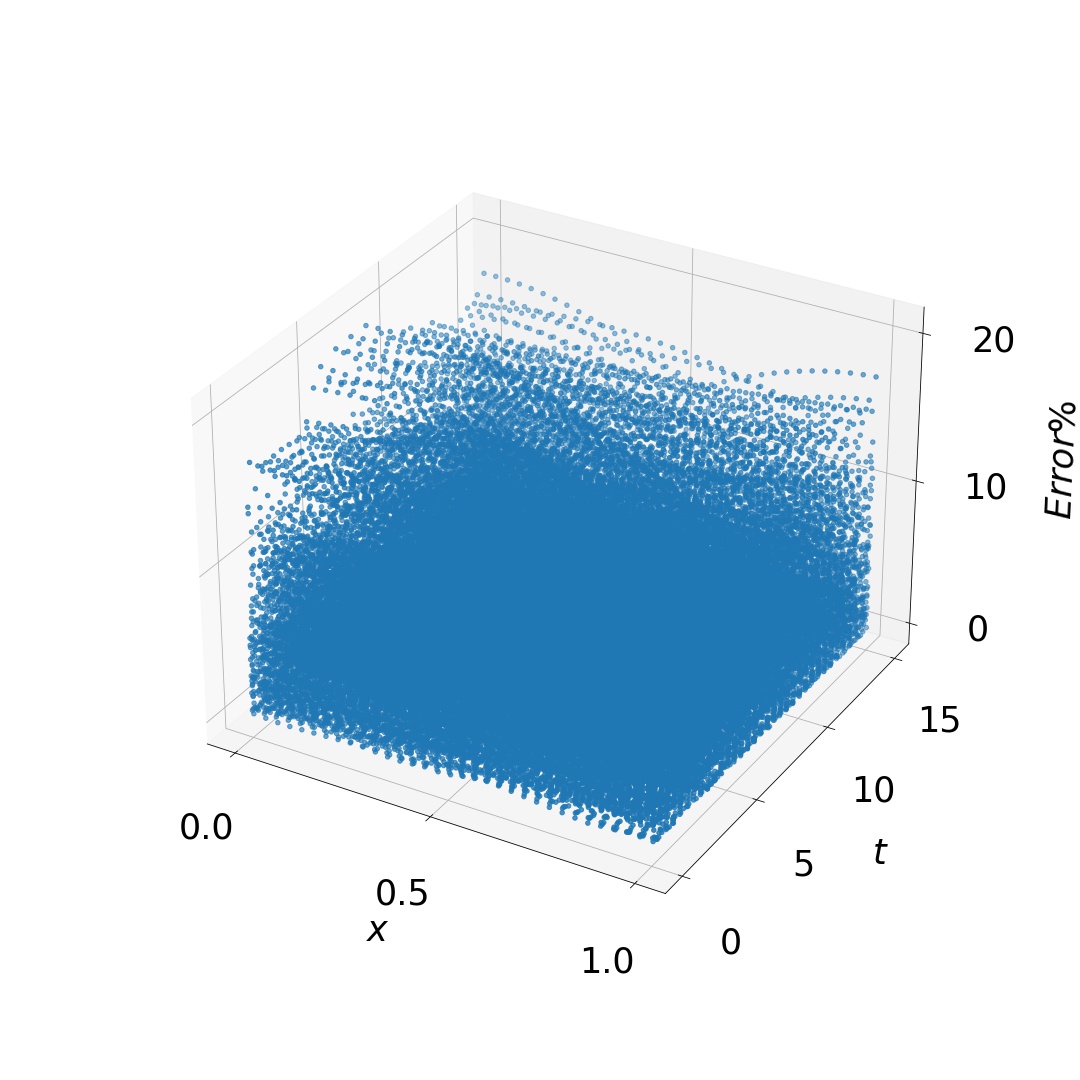}
		\caption{$\%$ Error with $i=5$ modes}
	\end{subfigure} \hspace{1em}
	\begin{subfigure}{0.4\linewidth}
		\includegraphics[width=\textwidth]{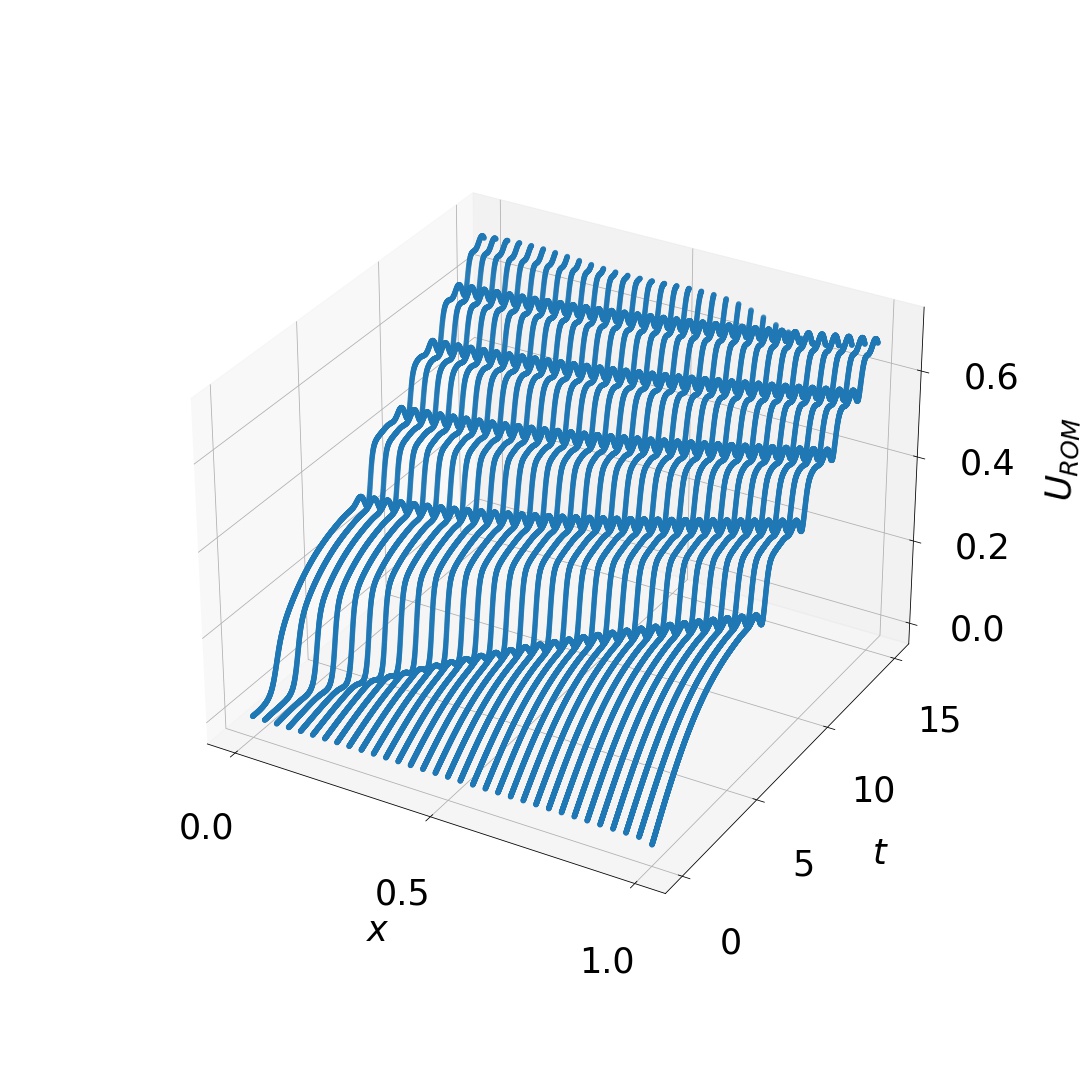}
		\caption{POD with $i=30$ modes. Model run time $t_{ROM}=0.227 s$}
	\end{subfigure}
	\begin{subfigure}{0.4\linewidth}
		\includegraphics[width=\textwidth]{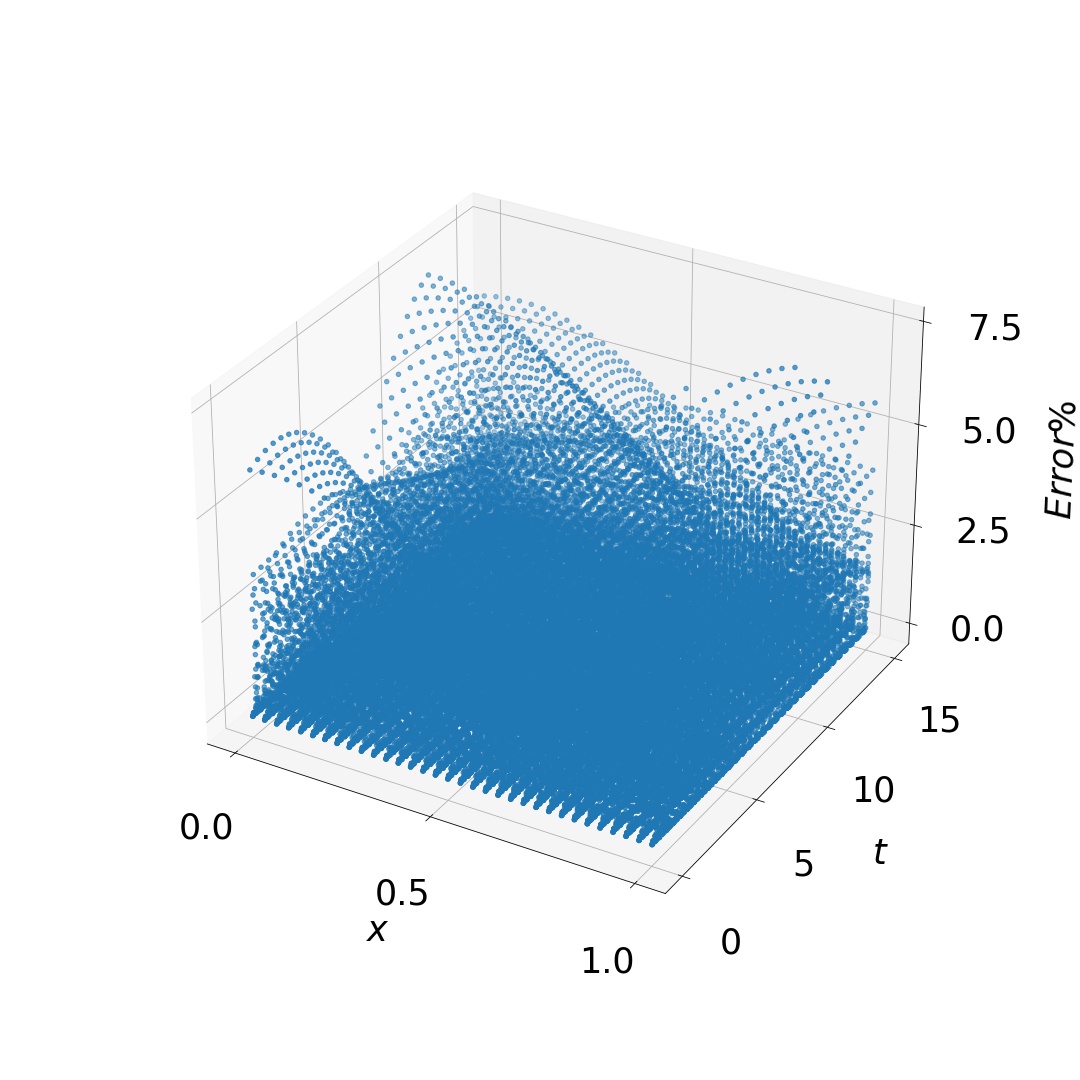}
		\caption{$\%$ Error with $i=30$ modes}
	\end{subfigure} \hspace{1em}
	\caption{Case 4 - Finite difference simulation with $2^{nd}$ order upwind scheme for convective term, longer simulation time and low diffusivity vs reduced order model via Galerkin POD.}
	\label{fig:gpod_case4}
\end{figure}

\begin{figure}[!htb]
	\centering
	\begin{subfigure}{0.4\linewidth}
		\includegraphics[width=\textwidth]{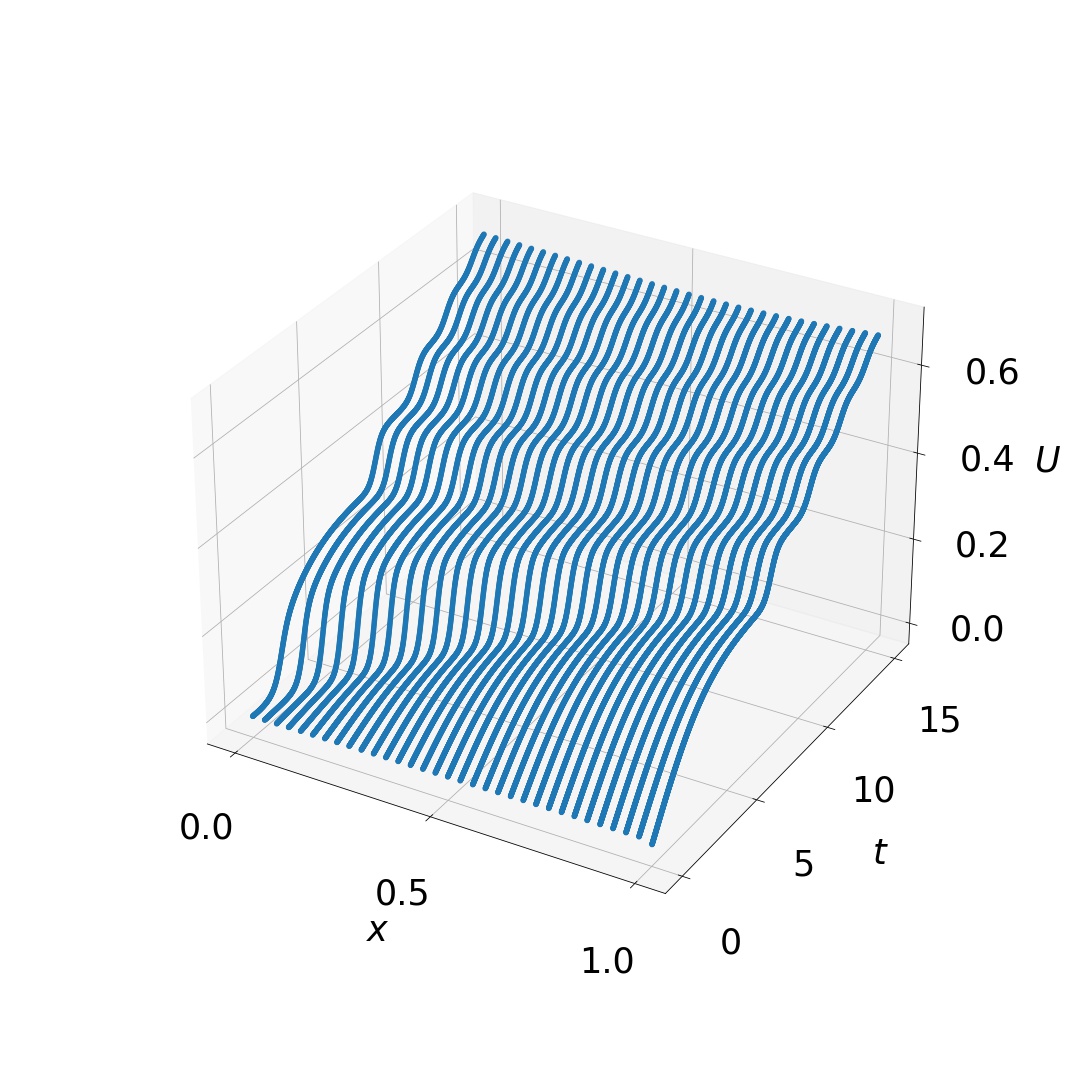} 
		\caption{Simulation with $t_{final}=15$, $N_x=32$, $CFL=0.2$, $\nu=0.0001$. Run time $t_{Sim}=0.485 s$}
		\label{fig_sim_case6}
	\end{subfigure} \hspace{1em}
	\begin{subfigure}{0.4\linewidth}
		\includegraphics[width=\textwidth]{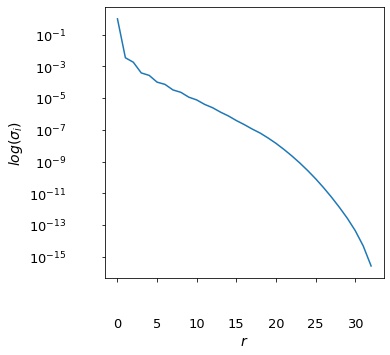} 
		\caption{Energy of dominant modes in the simulation}
		\label{fig_energy_case6}
	\end{subfigure} \\
	\begin{subfigure}{0.4\linewidth}
		\includegraphics[width=\textwidth]{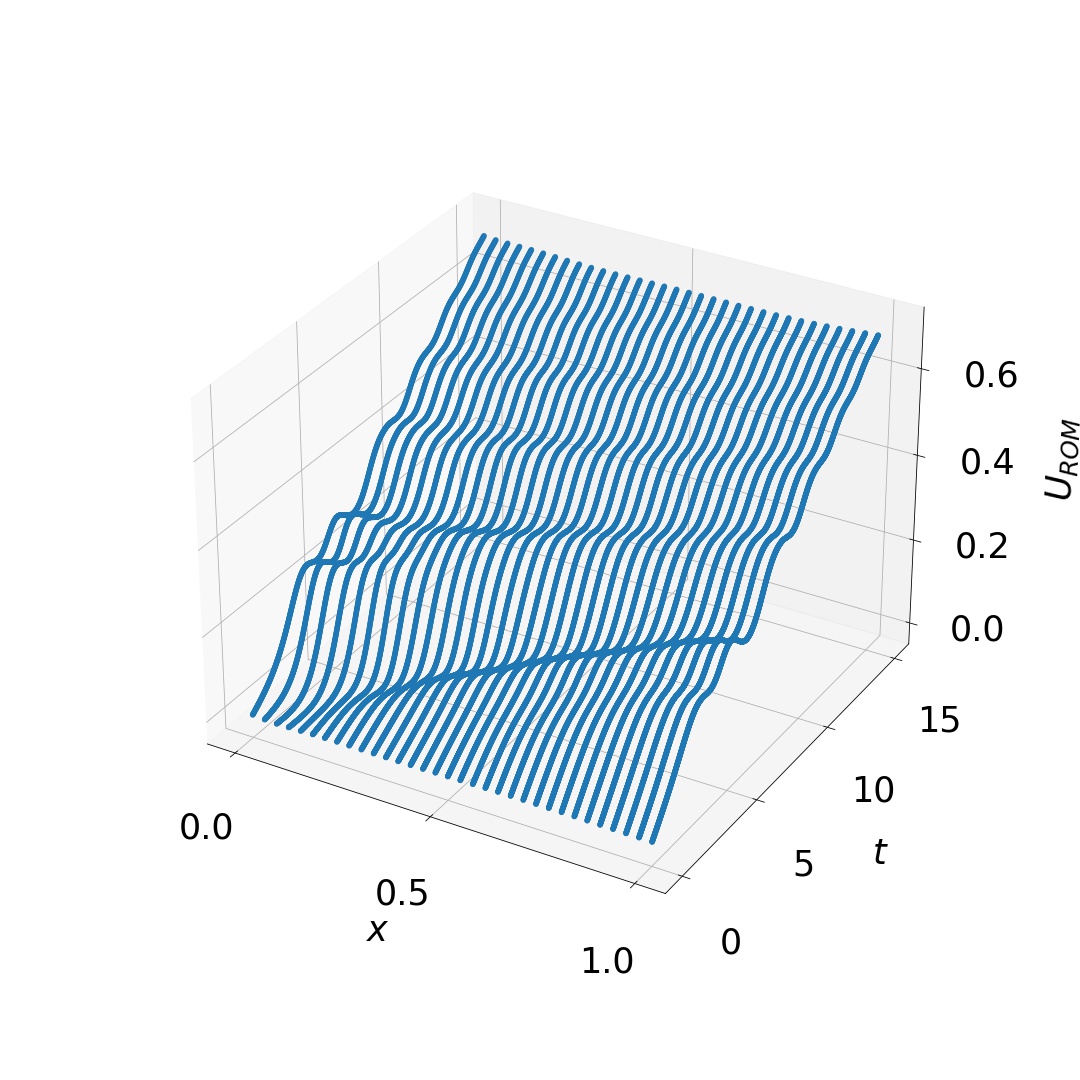}
		\caption{POD with $i=5$ modes. Model run time $t_{ROM}=0.06 s$}
	\end{subfigure} \hspace{1em}
	\begin{subfigure}{0.4\linewidth}
		\includegraphics[width=\textwidth]{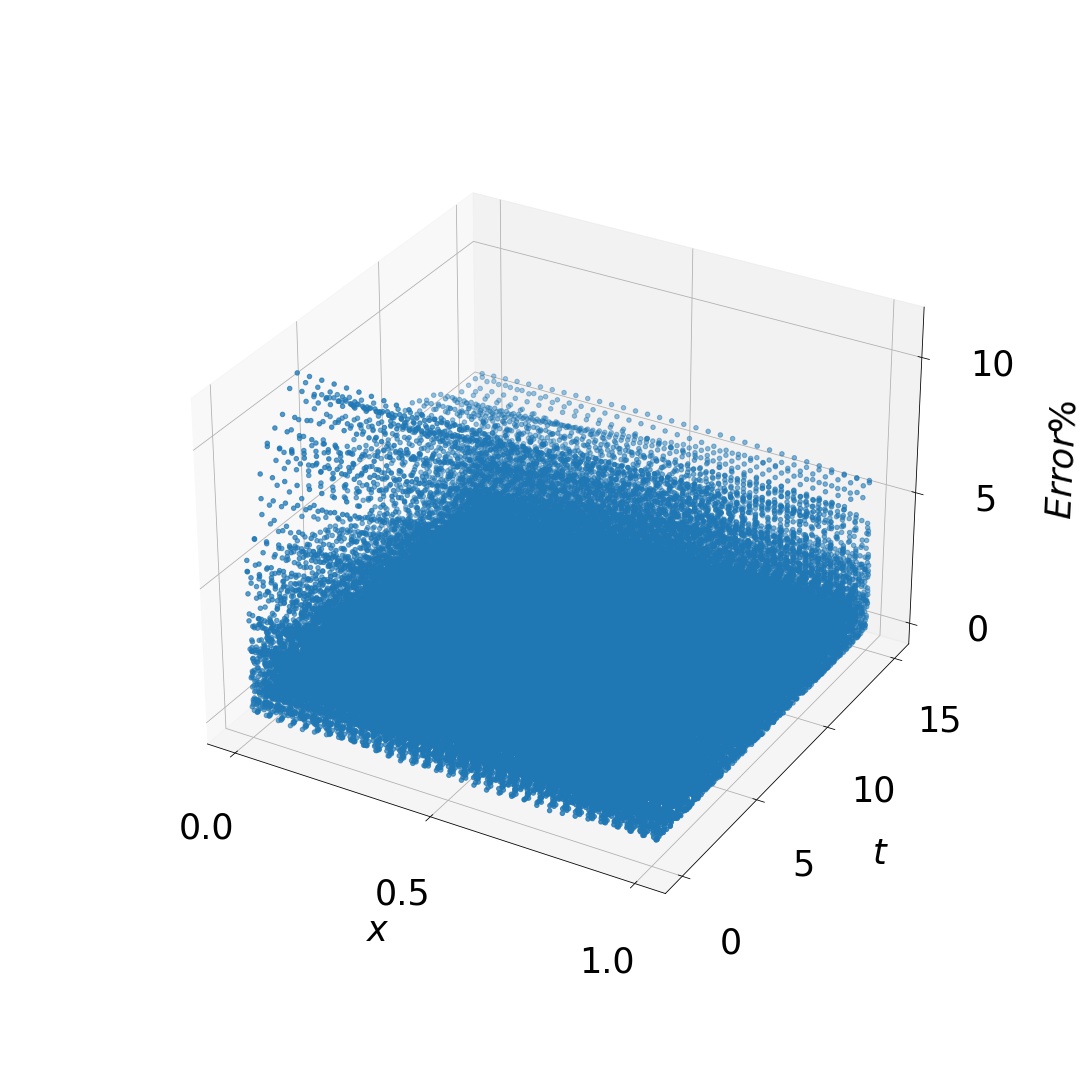}
		\caption{$\%$ Error with $i=5$ modes}
	\end{subfigure} \hspace{1em}
	\begin{subfigure}{0.4\linewidth}
		\includegraphics[width=\textwidth]{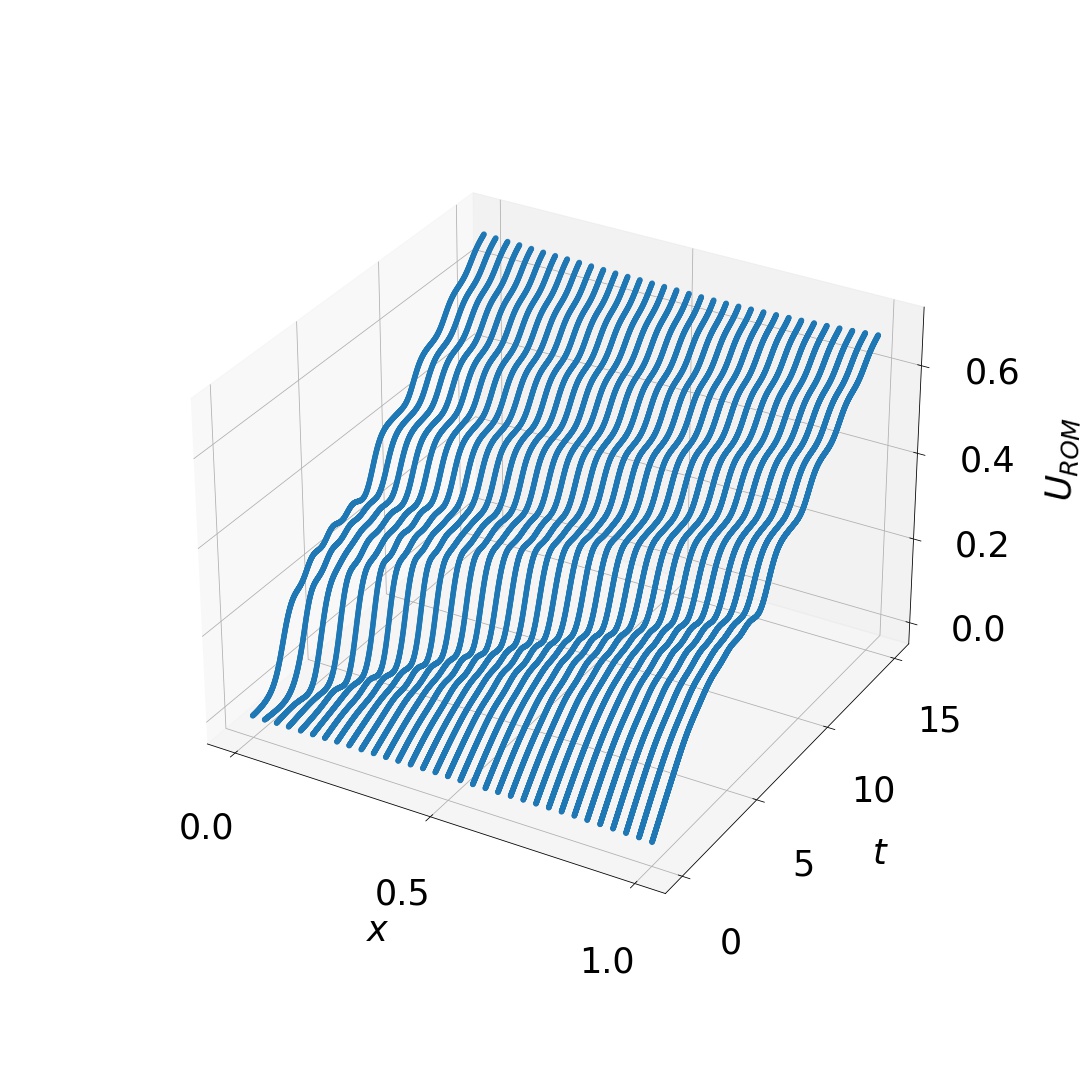}
		\caption{POD with $i=10$ modes. Model run time $t_{ROM}=0.093 s$}
	\end{subfigure}
	\begin{subfigure}{0.4\linewidth}
		\includegraphics[width=\textwidth]{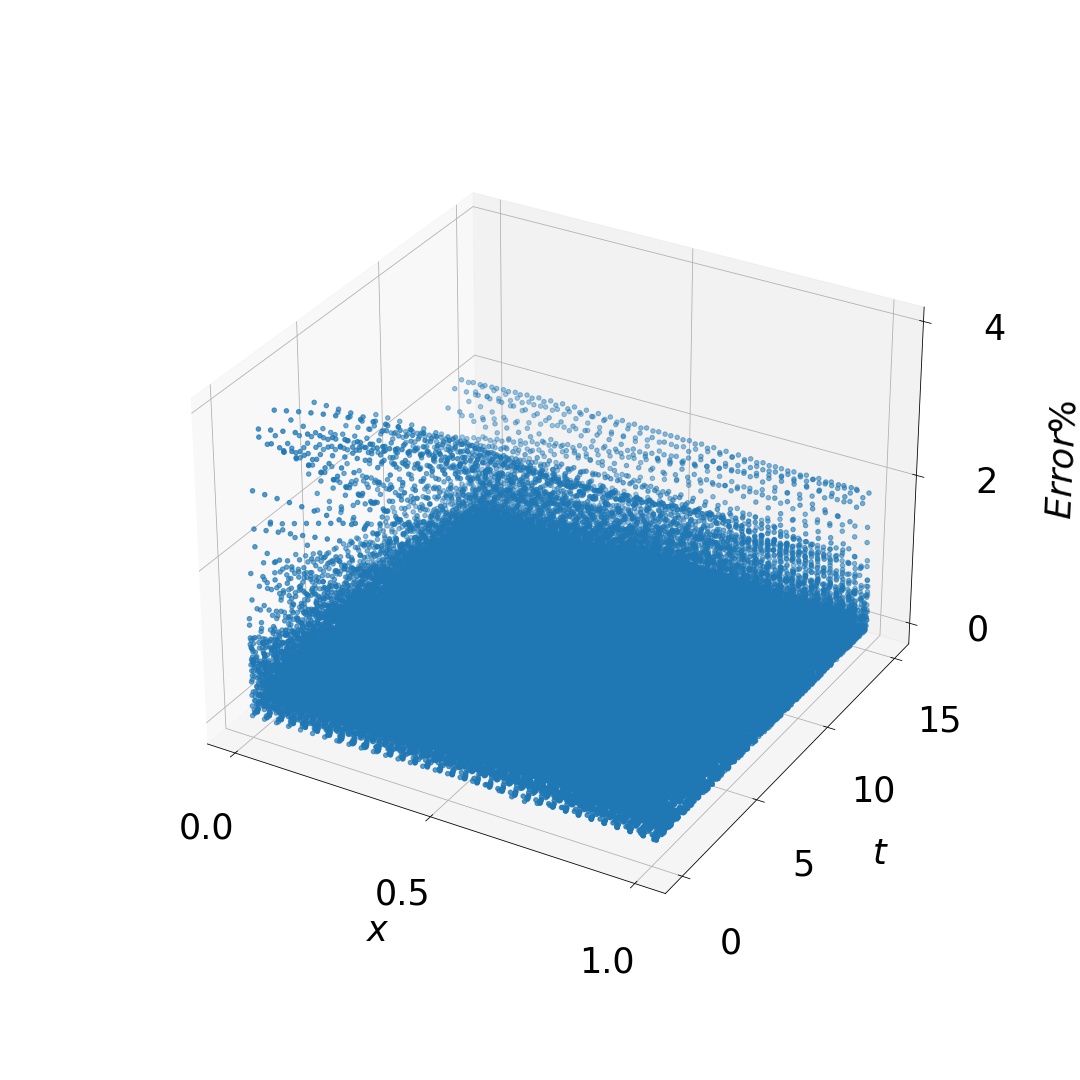}
		\caption{$\%$ Error with $i=10$ modes}
	\end{subfigure} \hspace{1em}
	\caption{Case 6 - Finite difference simulation with $1^{st}$ order upwind scheme for convective term, longer simulation run time and low diffusivity vs reduced order model via Galerkin POD.}
	\label{fig:gpod_case6}
\end{figure}

The results presented in this section evaluate the accuracy of POD-ROM for the following cases, Case 1: simulation with high diffusivity and short simulation run time ($t_{final}$), Case 2: simulation with high diffusivity and longer simulation run time, Case 3: simulation with low diffusivity and short run time, Case 4: simulation with low diffusivity and longer run time, and Case 5: heat equation with only diffusion term (Burgers' equation without convection term). For each case, the percentage error between the simulation and the model data is computed. For the case with high diffusivity and short simulation time (Fig.\ref{fig:gpod_case1}), it is observed that a satisfactory ROM can be created with as few as 1-3 dominant modes. However, with 5 dominant modes, the original data can be reconstructed with negligible error. At longer simulation times as in Case 2 (Fig. \ref{fig:gpod_case2}), a larger number of modes are required to reconstruct the data. For Case 1, it is observed that the POD-ROM constructed with 1 dominant mode between $t=0 s$ and $t=0.5 s$ captures the fluctuation in the field, but the POD-ROM constructed with 1 dominant mode in Case 2 does not capture the fluctuations between $t=0 s$ and $t=0.5 s$. This occurs because the POD modes are constructed from the data obtained from a specific simulation. As a result the POD modes and the distribution of energy across the modes vary for different data. At lower diffusivities and shorter simulation times as in Case 3 (Fig. \ref{fig:gpod_case3}), the trends similar to Case 1, where the data can be reconstructed with a few dominant modes, are observed. However, with increase in simulation time as in Case 4 (Fig. \ref{fig:gpod_case4}), large transient fluctuations are observed. At low diffusivity, Burger's equation is driven primarily by the convection term. If a parallel between this equation and the Navier Stokes equation is drawn, this is similar to a case at high Reynolds number. In such a scenario a large range of spatial and temporal scales exist, and the energy distribution is no longer dominated by a few modes. As a result, the POD-ROM constructed even with a large number of modes ($i= 5$ or $i= 30$) is inaccurate. The differences between the original data and ROM are observed to accentuate at simulation times $t > 10$. For 1D heat equation (Burgers' equation without the convection term), it is observed that POD-ROM can be constructed with only 1-3 modes even when the diffusivity is very small. Convection term is the primary source of non-linearity and fluctuations in Burgers' equation and POD is a great tool for construction of ROM for linear PDEs. An additional case (Case 6 - Fig.\ref{fig:gpod_case6}) is evaluated where a highly diffusive first order upwind scheme is used for discretization of the convection term. In this case, the fluctuations in the original data are observed to dissipate at longer simulation times. As a result, the energy distribution has a few dominant modes and the POD-ROM reconstructed with the limited number of modes is still accurate. This however is not an elegant solution since the original data from the simulation is highly diffusive and unreliable. 
\begin{table}[htb]
\centering
\begin{adjustbox}{width=0.6\textwidth}
\begin{tabular}{ccc} \hline
Case       & \multicolumn{2}{c}{$t_{ratio}=t_{ROM} / t_{Sim}$} \\ \hline
\multicolumn{3}{c}{}                                    \\
Case 1 & 0.76 $(i=1)$                & 1.03 $(i=5)$               \\
Case 2 & 0.076 $(i=1)$               & 0.15 $(i=5)$               \\
Case 3 & 0.57 $(i=1)$                & 1.12 $(i=5)$               \\
Case 4 & 0.31 $(i=20)$               & 0.42 $(i=33)$              \\
Case 6 & 0.12 $(i=5)$                & 0.2 $(i=10)$               \\ \hline
\end{tabular}
\end{adjustbox}
\caption{Computation time required for $l$ modes using POD method. $t_{ROM}$ represents the time required to run POD-ROM, $t_{Sim}$ represents the time required to run the finite difference simulation.}
\label{tab:comptime}
\end{table}
\section{Conclusions}
This work explores the development of ROM for 1D Burgers' equation using POD method. POD is very useful in creating fast and accurate ROM for a wide range of linear and non-linear equations. The primary idea of this work is to provide a feel for the method and present the material in an accessible format. It is the hope of the author that the reader is able to develop an appreciation for the approach and is able to apply it to more complex systems for reduction of computational times. 

\bibliographystyle{unsrt}
\bibliography{Ref}

\begin{thebibliography}{10}

\bibitem{burgers1}
Mikel Landajuela.
\newblock Burgers equation.
\newblock {\em BCAM Internship report: Basque Center for Applied Mathematics},
  2011.

\bibitem{burgers2}
Jusoh Yacob, Shahirah Zulkifli, Nur Zainuddin, and Siti Rosly.
\newblock Mathematical modelling of burger's equation applied in traffic flow.
\newblock Technical report, Kolej Teknologi Darulnaim, 10 2017.

\bibitem{burgers3}
A~Salih.
\newblock Burgers’ equation.
\newblock {\em Indian Institute of Space Science and Technology,
  Thiruvananthapuram}, 2016.

\bibitem{dissertation}
Neelakantan Padmanabhan.
\newblock {\em On High Pressure Real Gas Turbulent Mixing Jets}.
\newblock PhD thesis, All Dissertations. 1981., 2017.

\bibitem{POD2}
Farshid Abbasi and Javad Mohammadpour.
\newblock Nonlinear model order reduction of burgers' equation using proper
  orthogonal decomposition.
\newblock In {\em 2015 American Control Conference (ACC)}, pages 583--588.
  IEEE, 2015.

\bibitem{POD3}
Fariduddin Behzad, Brian~T Helenbrook, and Goodarz Ahmadi.
\newblock On the sensitivity and accuracy of
  proper-orthogonal-decomposition-based reduced order models for burgers
  equation.
\newblock {\em Computers \& Fluids}, 106:19--32, 2015.

\bibitem{POD5}
Eli Shlizerman, Edwin Ding, Matthew~O Williams, and J~Nathan Kutz.
\newblock The proper orthogonal decomposition for dimensionality reduction in
  mode-locked lasers and optical systems.
\newblock {\em International Journal of Optics}, 2012, 2012.

\bibitem{POD6}
Karl Kunisch and Stefan Volkwein.
\newblock Control of the burgers equation by a reduced-order approach using
  proper orthogonal decomposition.
\newblock {\em Journal of optimization theory and applications},
  102(2):345--371, 1999.

\bibitem{POD7}
Zhendong Luo, Xiaozhong Yang, and Yanjie Zhou.
\newblock A reduced finite difference scheme based on singular value
  decomposition and proper orthogonal decomposition for burgers equation.
\newblock {\em Journal of Computational and Applied Mathematics},
  229(1):97--107, 2009.

\bibitem{POD8}
Alfio Quarteroni, Gianluigi Rozza, et~al.
\newblock {\em Reduced order methods for modeling and computational reduction},
  volume~9.
\newblock Springer, 2014.

\bibitem{POD9}
Julien Weiss.
\newblock A tutorial on the proper orthogonal decomposition.
\newblock In {\em AIAA Aviation 2019 Forum}, page 3333, 2019.

\bibitem{POD10}
GRAU Katrin.
\newblock Applications of the proper orthogonal decomposition method.
\newblock Technical report, WN/CFD/07/97, 1997.

\bibitem{PCA1}
Jonathon Shlens.
\newblock A tutorial on principal component analysis.
\newblock {\em arXiv preprint arXiv:1404.1100}, 2014.

\bibitem{POD1}
Omer San and Traian Iliescu.
\newblock Proper orthogonal decomposition closure models for fluid flows:
  Burgers equation.
\newblock {\em arXiv preprint arXiv:1308.3276}, 2013.

\bibitem{POD4}
Jeff Borggaard, Zhu Wang, and Lizette Zietsman.
\newblock A goal-oriented reduced-order modeling approach for nonlinear
  systems.
\newblock {\em Computers \& Mathematics with Applications}, 71(11):2155--2169,
  2016.

\bibitem{POD_DEIM1}
Saifon Chaturantabut and Danny~C Sorensen.
\newblock Discrete empirical interpolation for nonlinear model reduction.
\newblock In {\em Proceedings of the 48h IEEE Conference on Decision and
  Control (CDC) held jointly with 2009 28th Chinese Control Conference}, pages
  4316--4321. IEEE, 2009.

\bibitem{POD_DEIM2}
Saifon Chaturantabut and Danny~C Sorensen.
\newblock Nonlinear model reduction via discrete empirical interpolation.
\newblock {\em SIAM Journal on Scientific Computing}, 32(5):2737--2764, 2010.

\bibitem{3blue1brown}
Grant Sanderson.
\newblock Linear algebra, an introduction to visualizing what matrices are
  really doing.
\newblock \url{https://www.3blue1brown.com/topics/linear-algebra}.
\newblock Accessed: 2023-03.

\bibitem{pacovideos}
Francisco Rodríguez~Fortuño.
\newblock Function vector spaces.
\newblock \url{https://www.youtube.com/watch?v=NvEZol2Q8rs&t=18s}.
\newblock Accessed: 2023-03.

\end{thebibliography}

\appendix
\section{Background in linear algebra and relevant terminologies} \label{linalgterm}
(Excerpts of lectures on the topic of linear algebra \cite{3blue1brown,pacovideos})
\begin{itemize}
    \item Vector space: A non-empty set $V$, of objects called vectors on which two operations (linear combinations) can be performed. Vector addition $(u+v \in V)$ and scalar multiplication $(\lambda u \in V)$.
    \item Span: A set of all possible vectors that can be formed by a linear combination of a set of vectors and scalars.
    \item Space: For a given vector space, if the tails of the vectors are at the origin and the tip of the vectors point at different locations, a grid space can be created such that each grid point is located at the tip of a given vector.
    \item Linear dependency: When two vectors point in the same coordinate direction and their span is restricted to a line or a plane, the vectors are called linearly dependent.
    \item Basis: A set of linearly independent vectors that span the full space. An infinite set of basis vectors exist, however physically only three vectors are perceived. For example, in a Cartesian coordinate the coordinate directions $x,y,z$ are defined as the basis.
    \item Linear transformation / Matrix-Vector multiplication: An operation, which when applied to a vector, it either rotates, scales, or performs a combination of the two to move it to a new location in space. A linear transformation is performed via matrix-vector multiplication $(Ax=b)$, where the matrix $A$ is the linear transformation applied to vector $x$ to form a new vector $b$. For a transformation to be defined as linear, the grids formed by the vector tips must remain parallel, evenly spaced and the origin must not move.
    \item Matrix multiplication: A linear transformation that is defined as a combination of two linear transformations. 
    \item Change of Basis: A linear transformation that translates a vector representation in one basis to a representation of the same vector in a new basis. 
    \item Eigenvector and Eigenvalue: When a linear transformation A changes the basis, some of the vectors $v$ in the space, remain on their span as an effect of the transformation. These vectors only scale in magnitude. This is defined as $Av=\lambda v$, where $v$ is defined as the Eigenvector and $\lambda$ is defined as Eigenvalue, which is a scaling factor. 
	\item Eigen decomposition: A factorization of a matrix $(A=V \Lambda V^{-1})$, where the matrix $A$ is represented a product of its Eigenvalues $\Lambda$ and Eigenvectors $V$. Eigen decomposition can only be applied to square matrices. 
	\item Singular value decomposition: A factorization of a real or complex matrix $(A=U \Sigma V^T)$ that is a generalization of Eigen decomposition but is applicable to non-square matrices as well. These factors may be physically interpreted as $U$ (orthogonal matrix - rotation), $\Sigma$ (diagonal matrix - stretching), and $V^T$ (orthogonal matrix - rotation)
	\item Function spaces \label{fnspace}: A set of all mathematical functions that have the same properties as vectors and vector spaces. For example, a function $f(x)$ that can be expressed a sum of two other functions $f(x)=u(x)+v(x)$ and a scalar product of another function $f(x)=\lambda w(x)$. From an abstract point of view, these functions are mathematically similar to a vectors and linear combinations can be applied to them to determine distance between them or to project them on one another. The difference in the analogy appears in the physical interpretation of the number of dimensions. It can be shown that an infinite number of linearly independent functions can be built from an infinite number of vectors to span the space. Therefore, functions exist in an infinite dimensional vector space. For example, a function $f(x)$ specified in the interval $[a,b]$, sampled at intervals $x_1,x_2,\ldots,x_N$ with corresponding function values of $f_1,f_2,\ldots,f_N$, can be expressed as, $$f(x)=\begin{bmatrix} f_1 \\
    \vdots\\
    f_N \end{bmatrix} =f_1 e_1+f_2 e_2+\dots+f_N e_N,$$
where, $e_i$ represents the basis of the function that takes the value of $1$ at the corresponding sampling point and $0$ elsewhere. A Dirac delta function may be used as a basis function. This product gives an approximation of the original function, and the approximation gets better as the number of samples are increased. In the limit $N \to \infty$, the basis functions become infinitely thin and infinitely many. Thus, the function over the interval is defined as,
$$f(x)=\int_{a} ^{b} f(\lambda) e_{\lambda} (x) d\lambda.$$
It must be noted that rectangular function is only one possible set of basis function. Other basis functions may also be used here to express the function, however not all basis functions exhibit the property of orthogonality. Using this definition of function space, a number of useful properties may be obtained such as projection of a function onto subspace span of another set of functions and determination of orthogonality of two functions. For an orthonormal basis ${e_1,e_2,\ldots,e_N }$,
$$a = \sum_i a_i e_i.$$
Inner product of two vectors: $\langle a,b \rangle=a_1 b_1^*+a_2 b_2^*+\ldots+a_N b_N^*$. Projection of a vector into a subspace of an orthonormal basis $W=span[w_1,w_2,\ldots,w_M]$,
$$w=\langle a,w_1 \rangle w_1+\langle a,w_2 \rangle w_2+\ldots+\langle a,w_M \rangle w_M.$$
\end{itemize}
\end{document}